\documentclass{amsart}
\usepackage{amssymb}
\input xy
\xyoption{all}

\newtheorem{theorem}{Theorem}[section]
\newtheorem{lemma}[theorem]{Lemma}
\newtheorem{proposition}[theorem]{Proposition}

\theoremstyle{definition}
\newtheorem{definition}[theorem]{Definition}
\newtheorem{eg}[theorem]{Example}

\newtheorem{remark}[theorem]{Remark}
\newtheorem{corollary}[theorem]{Corollary}

\numberwithin{equation}{section}

\begin{document}

\title{Characterizing Liminal And Type I Graph $C^*$-Algebras}

\author{Menassie Ephrem}
\address{Department of Mathematics and Statistics, Arizona State
University, Tempe, AZ 85287-1804}
\email{menassie@asu.edu}





\date{May, 2002.}



\begin{abstract}

We prove that the $C^*$-algebra of a directed graph $E$ is liminal
iff the graph satisfies the finiteness condition: if $p$ is an
infinite path or a path ending with a sink or an infinite emitter,
and if $v$ is any vertex, then there are only finitely many paths
starting with $v$ and ending with a vertex in $p$. Moreover,
$C^*(E)$ is Type I precisely when the circuits of $E$ are either
terminal or transitory, i.e., $E$ has no vertex which is on
multiple circuits, and $E$ satisfies the weaker condition: for any
infinite path $\lambda$, there are only finitely many vertices of
$\lambda$ that get back to $\lambda$ in an infinite number of
ways.
\end{abstract}

\maketitle



\section{Introduction}

A directed graph $E = (E^0, E^1,~~ o,~~ t)$ consists of a
countable set $E^0$ of vertices and $E^1$ of edges, and maps $o,t:
E^1 \rightarrow E^0$ identifying the origin (source) and the
terminus (range) of each edge. The graph is row-finite if each
vertex emits at most finitely many edges. A vertex is a sink if it
is not an origin of any edge.   A vertex $v$ is called singular if
it is either a sink or emits infinitely many edges.  A path is a
sequence of edges $e_1e_2\ldots e_n$ with $t(e_i) = o(e_{i+1})$
for each $i = 1,2,\ldots , n-1$. An infinite path is a sequence
$e_1e_2\ldots$ of edges with $t(e_i)=o(e_{i+1})$ for each $i$.

For a finite path $p=e_1e_2 \ldots e_n$, we define $o(p):=o(e_1)$
and $t(p):= t(e_n)$. For an infinite path $p=e_1e_2 \ldots$, we
define $o(p):=o(e_1)$.  We regard vertices as paths of length
zero, and hence if $v \in E^0, o(v)=v=t(v).$

$E^*=\bigcup_{n=0}^\infty E^n,$ where $E^n$ := \{$p:p$ is a path
of length n\}.

$E^{**}~~ :=~~ E^* \cup E^\infty$, where $E^\infty$ is the set of
infinite paths.

A Cuntz-Krieger $E$-family consists of mutually orthogonal
projections $\{ p_v : v \in E^0\}$ and partial isometries $\{s_e :
e \in E^1\}$ satisfying:

\begin{enumerate}

\item $p_{t(e)}=s_e^*s_e~~ \forall e \in E^1$.\\
\\

\item $\displaystyle{\sum_{e \in F}s_es_e^* \leq
p_v~~ \forall v \in E^0}$ and for any finite subset $F$ of $\{e
\in
E^1: o(e) = v\}$.\\
\\

\item $\displaystyle{\sum_{o(e)=v}s_es_e^* = p_v}$ for each
non-singular vertex $v \in E^0$.\\
\end{enumerate}

The graph $C^*$-algebra $C^*(E)$ is the universal $C^*$-algebra
generated by a Cuntz-Krieger $E$-family $\{s_e, p_v\}$.

For a finite path $\mu = e_1e_2 \ldots e_n$, we write $s_\mu$ for
$s_{e_1}s_{e_2}\ldots s_{e_n}$.

Since the family $\{s_\mu s_\nu^*:\mu, \nu \in E^*\}$ is closed
under multiplication, we have:

$$C^*(E) = \overline{span}\{s_\mu s_\nu^*:\mu, \nu \in E^*~~ and~~
t(\mu) = t(\nu)\}$$

The outline of the paper is as follows.  In section 2 we introduce
the basic notations and conventions we will use throughout the
paper. Section 3 deals with row-finite graphs with no sinks.  We
begin the section by defining a property of a graph which we later
prove to characterize liminal graph $C^*$-algebras when the graph
has no singular vertices. Section 4 provides us with a proposition
that gives us a method on how to obtain the largest liminal ideal
of a $C^*$-algebra of a row-finite graph with no sinks.  In
sections 5, respectively 6 , with the use of `desingularizing
graphs' of \cite{Dr}, we generalize the results of sections 3,
respectively 4, to arbitrary graphs. In section 7 we give a
characterization for type I graph $C^*$-algebras.  We finish the
section with a proposition on how to obtain the largest type I
ideal of a graph $C^*$-algebra.

I would like to thank my advisor Dr. Jack Spielberg.  He
consistently and tirelessly inspired me to explore the subject.
His help has been indispensable.

\section{Preliminaries}

Given a directed graph $E$, we write $v \geq w$ if there is a
directed path from $v$ to $w$.

For a directed graph $E$, we say $H \subseteq E^0$ is hereditary
if $v \in H$ and $v \geq w$ imply that $w \in H$. We say $H$ is
saturated if $v$ is not singular and $\{w \in E^0:v \geq w\}
\subseteq H$ imply that $v \in H$.

If $z \in \mathbb{T}$, then the family $\{zs_e, p_v\}$ is another
Cuntz-Krieger $E$-family with which generates $C^*(E)$, and the
universal property gives a homomorphism $\gamma_z:C^*(E)
\rightarrow C^*(E)$ such that $\gamma_z(s_e) = zs_e$ and
$\gamma_z(p_v) = p_v$.  $\gamma$ is a strongly continuous action,
called \textit{gauge action}, on $C^*(E)$.  See \cite{B} for
details.

Let $E$ be a row-finite directed graph, let $I$ be an ideal of
$C^*(E)$, and let $H=\{v:p_v \in I\}$.  In [\cite{B} Lemma 4.2]
they proved that $H$ is a hereditary saturated subset of $E^0$.
Moreover, if $I_H := \overline{span} \{S_\alpha S_\beta^*: \alpha,
\beta \in E^*~~ and~~ t(\alpha)= t(\beta) \in H\}$, the map $H
\longmapsto I_H$ is an isomorphism of the lattice of saturated
hereditary subsets of $E^0$ onto the lattice of closed
gauge-invariant ideals of $C^*(E)$ [\cite{B} Theorem 4.1 (a)].
Letting $F:= F(E \setminus H) =$ the sub-graph of $E$ that is
gotten by removing $H$ and all edges that point into $H$, it is
proven in [\cite{B} Theorem 4.1(b)] that $C^*(F) \cong
C^*(E)/I_H$. In the event that $I$ is not a gauge-invariant ideal,
we only get $I_H \subsetneqq I$.

\bigskip

\noindent We will use the following notations and conventions:

\begin{itemize}

\item[-] Every path we take is a directed path.

\item[-] A circuit in a graph $E$ is a finite path $p$ with $o(p) =
t(p)$.  We save the term loop for a circuit of length 1.

\item[-] We say that a circuit is terminal if it has no exits, and
a circuit is transitory if it has an exit and no exit of the
circuit gets back to the circuit.

\item[-] $\Lambda_E~~ :=~~ \{v \in E^0: v~~ is~~ a~~ singular~~ vertex\}$

\item[-] $\Lambda_E^*~~ :=~~ t^{-1}(\Lambda_E) \cap E^*$ i.e.,
$\Lambda_E^*$ is the set of paths ending with a singular vertex.
When there are no ambiguities, we will just use $\Lambda^*$.

\item[-] We say $v$ gets to $w$ (or reaches $w$) if there is a path from
$v$ to $w$.

\item[-] We say $v$ gets to a path $p$ if $v$ gets to a vertex in $p$.

\item[-] For a subset $S$ of $E^0$, we write $S \geq v$ if $w \geq v~~
\forall w \in S$.

\item[-] For a subset $H$ of $E^0$, we write $Graph(H)$ to refer to the
sub-graph of $E$ whose set of vertices is $H$ and whose edges are
those edges of $E$ that begin and end in $H$.

\item[-] $V(v):= \{w \in E^0: v \geq w\}$.

\item[-] $E(v) := Graph(V(v))$.  i.e., $E(v)$ is that part of the sub-graph
of $E$ that the vertex $v$ can `see'.  Accordingly we use $F(v)$,
etc. when the graph is $F$, etc.

\item[-] For $v \in E^0$ let $\Delta(v) := \{e \in E^1 : o(e) = v\}$.

\item[-] For a finite subset $F$ of $\Delta(v)$, we write $\displaystyle
V(v; F) := \{v\} \cup \bigcup_{e \in \Delta(v) \setminus F}
V(t(e))$.

\item[-] $E(v; F):= Graph(V(v; F))$.

\item[-] For a hereditary subset $H$ of $E^0$, we write $\overline{H}$ to
refer to the saturation of $H$ = the smallest saturated set
containing $H$.  Notice that $\overline{H}$ is hereditary and
saturated.

\item[-] For any path $\lambda$, $\lambda^0$ will denote the vertices of
$\lambda$.

\item[-] As was used above, $F(E \setminus H)$ will denote the sub-graph of
$E$ that is gotten by removing $H$ and all edges that point into
$H$

\item[-] We use $\mathcal{K}$ to denote the space of compact operators
on an (unspecified) separable Hilbert space.
\end{itemize}

\section{Liminal $C^*$-Algebras of Graphs with no singular vertices.}

We begin the section by a definition.

\begin{definition}\label{maxi.tail} A subset $\gamma$ of $E^0$ is
called a maximal tail if it satisfies the following three
conditions.

\begin{enumerate}

\item[(a)] for any $v_1, v_2 \in \gamma$ there exists $z \in \gamma$
s.t. $v_1 \geq z$ and $v_2 \geq z$ .

\item[(b)] for any $v \in \gamma~~ \exists e \in E^1~~ s.t.~~ o(e) =
v$ and $t(e) \in \gamma$.

\item[(c)] $v \geq w$ and $w \in \gamma$ imply that $v \in \gamma$ .
\end{enumerate}
\end{definition}

We will prove a result similar to (one direction of) [\cite{B}
Proposition 6.1] with a weaker assumption on the graph $E$ and a
weaker assumption on the ideal.

\begin{lemma}\label{lemm.1}  Let $E$ be a row-finite graph with no
sinks.  If $I$ is a primitive ideal of $C^*(E)$ and $H = \{v \in
E^0: p_v \in I\}$, then $\gamma=E^0 \setminus H$ is a maximal
tail.
\end{lemma}

\begin{proof}  By [\cite{B} Lemma 4.2] $H$ is hereditary saturated.
The complement of a hereditary set satisfies (c).  Since $E$ has
no sinks, and $H$ is saturated, $\gamma$ satisfies (b).  We prove
 (a).  Let $v_1,~~ v_2 \in \gamma$ and let $H_i = \{v \in \gamma:v_i
 \geq v\}$.  We will first show that $\overline{H_1} \cap
\overline{H_2} \neq \emptyset$.  Let
 $F=F(E \setminus H)$. For each $i$, $I_{\overline{H_i}}$ is a
non-zero ideal of $C^*(F)
 \cong C^*(E)/I_H$, hence is of the form $I_i/I_H$, and $p_{v_i} +
 I_H \in I_{\overline{H_i}}$. Since each $I_{\overline{H_i}}$ is
 gauge-invariant, so is $I_{\overline{H_1}}
 \cap I_{\overline{H_2}}$.  Therefore $I_{\overline{H_1}} \cap
 I_{\overline{H_2}} = I_{\overline{H_1} \cap \overline{H_2}}$.  If
$\overline{H_1} \cap \overline{H_2} = \emptyset$ then
$I_{\overline{H_1}} \cap I_{\overline{H_2}} = \{0\} \subseteq
I/I_H$.  But $I/I_H$ is a primitive ideal of $C^*(E)/I_H$
therefore $I_1/I_H \subseteq I/I_H$ or $I_2/I_H \subseteq I/I_H$.
WLOG let $I_1/I_H \subseteq I/I_H$ hence $p_{v_1} + I_H \in I/I_H$
implying that $p_{v_1} \in I_H$ or $p_{v_1} \in I \setminus I_H$.
But $p_{v_1} \in I \setminus I_H$ is a contradiction to the
construction of $H$, and $p_{v_1} \in I_H$, which implies that
$v_1 \in H$, which is again a contradiction to $v_1 \in \gamma =
E^0 \setminus H$. Therefore $\overline{H_1} \cap \overline{H_2}
\neq \emptyset$. Let $y \in \overline{H_1} \cap \overline{H_2}$.
Applying [\cite{B} Lemma 6.2] to $F$ and $v_1$ shows that $\exists
x \in E^0 \setminus H~~ \textrm{s.t.}~~ y \geq x$ and $v_1 \geq
x$. Since $y \in \overline{H_2}$ and $\overline{H_2}$ is
hereditary, $x \in \overline{H_2}$. Applying [\cite{B} Lemma 6.2]
to $F$ and $v_2$ shows that $\exists z \in E^0 \setminus H~~
\textrm{s.t.}~~ y \geq z$ and $v_2 \geq z$.  Thus $v_1 \geq z,$
and $v_2 \geq z$ as needed.
\end{proof}

Now, we prove that for a row-finite graph $E$ with no sinks,
$C^*(E)$ is liminal precisely when the following finiteness
condition is satisfied: for any vertex $v$ and any infinite path
$\lambda$, there is only a finite number of ways to get to
$\lambda$ from $v$.

To state the finiteness condition more precisely, we will use the
equivalence relation defined in [\cite{S} Definition 1.8].

If  $p = e_1e_2 \ldots$ and $q = f_1f_2 \cdots \in E^\infty$, we
say that $p \sim q$ iff $\exists j, k$ so that $e_{j+r} = f_{k+r}$
for $r \geq 0$. i.e., iff $p$ and $q$ (eventually) share the same
tail.

We use $[p]$ to denote the equivalence class containing $p$.

\begin{definition}\label{defn.1}
A row-finite directed graph $E$ that has no sinks is said to
satisfy condition $(M)$ if for any $v \in E^0$ and any $[p] \in
E^\infty/\sim$ there is only a finite number of representatives of
$[p]$ that begin with v.
\end{definition}

We note that $E$ satisfies condition $(M)$ implies that every
circuit in $E$ is terminal.

\begin{lemma}\label{lemm.2} Let $E$ be a row-finite directed graph
with no sinks that satisfies condition $(M)$.  Let $F$ be a
sub-graph of $E$ so that $F^0$ is a maximal tail. If $F$ has a
circuit, say $\alpha$, then the saturation of $\alpha^0$,
$\overline{\alpha^0}$, is equal to $F^0$.
\end{lemma}

\begin{proof}
Let $v_\alpha$ be a vertex of $\alpha$. Since $\alpha$ is
terminal, $v_\alpha \geq z$ implies that $z$ is in $\alpha^0$.
Also, for each $w \in F^0$, by (a) of Definition \ref{maxi.tail},
there exists $z \in F^0$ s.t. $w \geq z$ and $v_\alpha \geq z$,
but $z$ is in $\alpha^0$ which implies that $z \geq v_\alpha$
Therefore $w \geq v_\alpha$, i.e., each vertex in $F^0$ connects
to $v_\alpha$ (via a directed path).

Now, assuming the contrary, let $v_1 \notin \overline{\alpha^0}$.
If $v_1$ is in a circuit, say $\beta$.  Then, by the previous
paragraph, $v_1 \geq v_\alpha$ hence either $\beta = \alpha$ or
$\beta$ has an exit.   But $v_1 \notin \overline{\alpha^0}$,
therefore $\beta = \alpha$ is not possible, and since $F$
satisfies condition $(M)$, $\beta$ can not have an exit. Thus
$v_1$ is not in a circuit.  Therefore $\exists e_1 \in
F^1~~{\textrm s.t.}~~ o(e_1) = v_1,~~ t(e_1) \notin
\overline{\alpha^0}.$  Let $v_2 = t(e_1)$. Inductively, $\exists
e_n \in F^1~~ \textrm{s.t.}~~ v_n= o(e_n),~~ t(e_n) = v_{n+1}
\notin \overline{\alpha^0}$. Look at the infinite path $e_1e_2
\ldots$.

Notice that the $v_i$'s are distinct and each $v_i \geq v_\alpha$.
Therefore there are infinitely many ways to get to $\alpha$ from
$v_1$, i.e., there are infinitely many representatives of
[$\alpha$] that begin with $v_1$, which contradicts to the
assumption that $E$ satisfies condition $(M)$.  Therefore $F^0 =
\overline{\alpha^0}$
\end{proof}

\begin{lemma}\label{lemm.3} Let $E$ be a row-finite directed graph
with no sinks that satisfies condition $(M)$.  Let $F$ be a
sub-graph of $E$ so that $F^0$ is a maximal tail. If $F$ has no
circuits then $F$ has a hereditary infinite path, say $\lambda$,
s.t. $F^0 = \overline{\lambda^0}$.
\end{lemma}

\begin{proof} Since $F$ has no sinks, it must have an infinite
path, say $\lambda$. Let $v_\lambda$ be a vertex in $\lambda$.  By
condition $(M)$, there are only a finite number of infinite paths
that begin with $v_\lambda$ and share a tail with $\lambda$.  By
going far enough on $\lambda$, there exists $w \in \lambda^0$ s.t.
$v_\lambda \geq w$ and $[\lambda]$ has only one representative
that begins with $w$.  By re-selecting $v_\lambda$ (to be $w$, for
instance) we can assume that there is only one representative of
$[\lambda]$ that begins with $v_\lambda$.  We might, as well,
assume that $o(\lambda) = v_\lambda$.

We will now prove that $\lambda^0$ is hereditary.  Suppose $u \in
F^0$ s.t. $v_\lambda \geq u$ and $u \notin \lambda^0$. Since $F^0$
is a maximal tail and since $F$ has no circuits, by (b) of
Definition \ref{maxi.tail} we can choose $w_1 \in F^0$ s.t.
$v_\lambda \geq w_1$ and $v_\lambda \neq w_1$. By (a) of
Definition \ref{maxi.tail} there exists $z_1 \in F^0$ s.t. $u \geq
z_1$, and $w_1 \geq z_1$. If $z_1 \in \lambda^0$ then we have two
ways to get to $\lambda$ from $v_\lambda$ (through $u$ and through
$w_1$) which contradicts to the choice of $v_\lambda$, hence $z_1
\notin \lambda^0$

Let $w_2 \in \lambda^0$ (far enough) so that $w_2 \ngeq z_1$. If
such a choice was not possible, we would be able to get to $z_1$
and hence to any path that begins with $z_1$ from $v_\lambda$ in
an infinite number of ways, contradicting condition $(M)$.

Again since $F^0$ is a maximal tail, there exists $z_2 \in F^0$
s.t. $w_2 \geq z_2$ and $z_1 \geq z_2$.  Notice that there are (at
least) two ways to get to $z_2$ from $v_\lambda$.  By inductively
choosing a $w_n \in \lambda^0$ and a $z_n \in F^0$ s.t. $w_n \ngeq
z_{n-1}, w_n \geq z_n$ and $z_{n-1} \geq z_n$, there are at least
$n$ ways to get to $z_n$ from $v_\lambda$ (one through $w_n$ and
$n-1$ through $z_{n-1}$).

We now form an infinite path $\alpha$ that contains $z_1, z_2
\ldots$ as (some of) its vertices that we can reach to, from
$v_\lambda$, in an infinite number of ways, which is again a
contradiction. Hence no such $u$ can exist.  Thus $\lambda$ is
hereditary.

We will now prove that $F^0 = \overline{\lambda^0}$. Assuming the
contrary, let $v_1 \notin \overline{\lambda^0}$. Then $\exists e_1
\in F^1~~{\textrm s.t.}~~ o(e_1) = v_1,~~ t(e_1) \notin
\overline{\lambda^0}.$   Inductively, let $v_n = t(e_{n-1})$, then
$\exists e_n \in F^1~~ \textrm{s.t.}~~ v_n= o(e_n),~~ t(e_n) =
v_{n+1} \notin \overline{\lambda^0}$. Consider the infinite path
$e_1e_2 \ldots$.

Notice that since $F^0$ is a maximal tail, for each $v_i~~ \exists
x_i$ s.t. $v_i \geq x_i$ and $v_\lambda \geq x_i$. But $\lambda^0$
is hereditary hence $x_i \in \lambda^0$, implying that each $v_i$
reaches $\lambda$. Therefore there are infinitely many ways to get
to $\lambda$ from $v_1$, i.e., there are infinitely many
representatives of [$\lambda$] that begin with $v_1$ which
contradicts to condition $(M)$. Therefore $F^0 =
\overline{\lambda^0}$
\end{proof}

\begin{lemma}\label{lemm.4} Let $I_H$ be a primitive ideal of
$C^*(E)$, where $H$ is a hereditary saturated subset of $E^0$. Let
$F = F(E \setminus H)$.  Then $F$ has no circuits.
\end{lemma}

\begin{proof} Note that $F^0$ is a maximal tail. And $C^*(F) \cong
C^*(E)/I_H$. Since $I_H$ is a primitive ideal of $C^*(E),~~ \{0\}$
is a primitive ideal of $C^*(F)$.

Suppose that $F$ has a circuit, say $\alpha$.  By Lemma
\ref{lemm.2}, $F^0 = \overline{\alpha^0}$.  Hence $C^*(F) \cong
I_{\alpha^0}$ = the ideal of $C^*(F)$ generated by $\{\alpha^0\}$.
Since $\alpha$ has no exits (is hereditary), by [\cite{K}
Proposition 2.1] $I_{\alpha^0}$ is Morita equivalent to
$C^*(\alpha)$ which is Morita equivalent to $C(\mathbb{T})$.  But
$\{0\}$ is not a primitive ideal of $C(\mathbb{T})$ implying that
$\{0\}$ is not a primitive ideal of $C^*(F)$ which is a
contradiction. Hence $F$ has no circuits.
\end{proof}

Hidden in the proofs of Lemma \ref{lemm.3} and Lemma \ref{lemm.4}
we have proven a (less relevant) fact: If a directed graph $E$
with no singular vertices satisfies condition $(M)$ and $F^0$ is a
maximal tail then $F$ has (essentially) one infinite tail, i.e.,
$F^\infty/\sim$ contains a single element.

\begin{remark}\label{remk.1}  Let $E_1$ be a sub-graph of a
directed graph $E_2$.  Applying [\cite{S} Theorem 2.34] and
[\cite{S} Corollary 2.33], we observe that $C^*(E_1)$ is a
quotient of a $C^*$-subalgebra of $C^*(E_2)$. (Letting $S_2 =
E_2^0 \setminus \{v \in E_2^0:v \textrm{ is a singular vertex}
\}$)
\end{remark}

\begin{lemma}\label{lemm.5} Let $E$ be a directed graph.  Suppose
all the circuits of $E$ are transitory and suppose $\exists\lambda
\in E^\infty$ s.t. the number of vertices of $\lambda$ that emit
multiple edges that get back to $\lambda$ is infinite.  Then
$C^*(E)$ is not Type I.
\end{lemma}

\begin{proof}

Let $v_1 \in \lambda^0$ s.t. $v_1$ emits (at least) two edges that
get back to $\lambda$.

choose a path $\alpha_1^1 = e_1e_2 \ldots e_{n_1}$, s.t. $e_1$ is
not in $\lambda,~~ o(e_1) = v_1$ and $t(\alpha_1^1) = v_2 \in
\lambda^0$. If $t(e_{n_1}) = v_1$, i.e. $e_1e_2 \ldots e_{n_1}$ is
a circuit we extend $e_1e_2 \ldots e_{n_1}$ so that $v_2$ is
further along $\lambda$ than $v_1$ is.

We might again extend $\alpha^1_1$ along $\lambda$, if needed, and
assume that $v_2$ emits (at least) two edges that get back to
$\lambda$.

Let $\alpha_1^2$ be the path along $\lambda$ s.t. $o(\alpha^2_1) =
v_1$ and $t(\alpha_1^2) = v_2.$  Inductively, choose $\alpha^1_k =
e_1e_2 \ldots e_{n_k}$ s.t. $e_1$ is not in $\lambda$, $o(e_1) =
v_k,~~ t(\alpha^1_k) = v_{k+1} \in \lambda^0$, by extending
$\alpha^1_k$, if needed, we can assume that $v_{k+1}$ is further
along $\lambda$ than $v_k$ and emits multiple edges that get back
to $\lambda$.  Let $\alpha^2_k$ be the path along $\lambda$ s.t.
$o(\alpha^2_k) = v_k$ and $t(\alpha^2_k) = v_{k+1}$. Now look
at the following sub-graph of $E$, call it $F$.

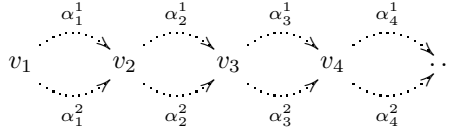
$$\xymatrix{v_1 \ar @{.>} @/^1pc/ [r]^{\alpha_1^1} \ar @{.>}
@/_1pc/[r]_{\alpha_1^2} & v_2 \ar @{.>} @/^1pc/[r]^{\alpha_2^1}
\ar @{.>} @/_1pc/[r]_{\alpha_2^2} & v_3 \ar @{.>}
@/^1pc/[r]^{\alpha_3^1} \ar @{.>} @/_1pc/[r]_{\alpha_3^2} & v_4
\ar @{.>} @/^1pc/[r]^{\alpha_4^1} \ar @{.>}
@/_1pc/[r]_{\alpha_4^2} & { } \ldots}$$

Now let $\{s_e, p_v: e \in F^1, v \in F^0\}$ be a Cuntz-Krieger
$F$-family.

Thus $C^*(F) = \overline{span}\{s_\mu s_\nu^*: \mu, \nu \in F^*~~
and~~ t(\mu) = t(\nu)\}.$

Let $F_k:= \overline{span}\{s_\mu s_\nu^*: \mu, \nu$ are paths
made up of $\alpha^r_i$'s (or just $v_k$) s.t. $t(\mu) = v_k =
t(\nu)\}$.

By [\cite{K} Corollary 2.3], $F_k \cong M_{N_k}(\mathbb{C})$ where
$N_k$ is the number of paths made up of $\alpha^r_i$'s (or just
$v_k$) ending with $v_k$, which is finite.

Also, if $s_\mu s_\nu^* \in F_k$ then $s_\mu s_\nu^* = s_\mu
p_{v_k} s_\nu^* = s_\mu s_{\alpha^1_k} s^*_{\alpha^1_k} s^*_\nu +
s_\mu s_{\alpha^2_k} s^*_{\alpha^2_k} s^*_\nu \in F_{k + 1}.$
Hence $F_k \subsetneqq F_{k+1}$.  Let
$\displaystyle{\mathcal{A}=\overline{\bigcup_{k=1}^\infty F_k}}$.
Then $\mathcal{A}$ is a $C^*$-subalgebra of $C^*(F)$. Since
$\mathcal{A}$ is a UHF algebra, it is not Type I. Therefore
$C^*(F)$ has a $C^*$-subalgebra that is not Type I and can not be
Type I. Since $F$ is a sub-graph of $E$, by Remark \ref{remk.1}
$C^*(E)$ has a $C^*$-subalgebra whose quotient is not Type I.
Therefore $C^*(E)$ is not Type I.
\end{proof}

It might be useful to keep following graph in mind when reading
Lemma \ref{lemm.51}, it can be viewed as a prototype of a graph
that satisfies the assumption of the lemma.

$$
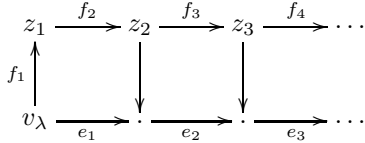
\xymatrix{z_1 \ar[r]^{f_2} & z_2 \ar[r]^{f_3} \ar[d] & z_3
\ar[r]^{f_4} \ar[d] &  {\ldots} \\
            v_\lambda \ar[r]_{e_1} \ar[u]^{f_1} & .
            \ar[r]_{e_2} & . \ar[r]_{e_3} & \ldots}$$

\begin{lemma}\label{lemm.51} Let $E$ be a directed graph, let
$\lambda = e_1e_2 \ldots$ be an infinite path in $E$, and let
$o(\lambda) = v_\lambda$.  Suppose:
\begin{enumerate}

\item $E$ has no circuits.

\item The number of representatives of $[\lambda]$ that begin with
$v_\lambda$ is infinite.

\item $v_\lambda$ is the only such vertex in $\lambda^0$.

\item $E = E(v_\lambda)$, and

\item  $\forall v \in E^0~~ \exists w \in \lambda^0~~
\textrm{s.t.}~~ v \geq w~~ (i.e., E^0 \geq \lambda^0)$.
\end{enumerate}

Then

 i) $\{0\}$ is a primitive ideal of $C^*(E)$, and

ii) $C^*(E)$ is not simple.
\end{lemma}

\begin{proof}
We prove i). First note that $E$ satisfies condition $(K)$ of
\cite{B}: every vertex lies on either no circuits or at least two
circuits.  This is because $E$ has no circuits. We will show that
$E^0$ is a maximal tail.  Since $E$ has no sinks, $E^0$ satisfies
(b) of Definition \ref{maxi.tail}, and clearly $E^0$ satisfies
(c). We will show that $E$ satisfies (a). Let $v_1, v_2 \in E^0$.
By (4) above, $\exists w_1, w_2 \in \lambda^0$ s.t. $v_i \geq
w_i$. Since $\lambda$ is is an infinite path, either $w_1 \geq
w_2$ or $w_2 \geq w_1$. WLOG let $w_2 \geq w_1$. We have $v_1 \geq
w_1$ and $v_2 \geq w_2 \geq w_1$, hence (a) is satisfied.
Therefore $E^0$ is a maximal tail and, by [\cite{B} Proposition
6.1], $I_\emptyset = \{0\}$ is a primitive ideal of $C^*(E)$.

We will prove ii). Since $E$ is row finite and since $v_\lambda$
gets to $\lambda$ infinitely often, $\exists f_1 \in E^1$ s.t.
$o(f_1) = v_\lambda$ and $z_1 := t(f_1)$ gets to $\lambda$
infinitely often.  Moreover there is no vertex in $\lambda^0$ that
gets to $\lambda$ infinitely often except $v_\lambda$ and $E$ has
no circuits, therefore $z_1 \notin \lambda^0$. Inductively,
$\exists f_{n+1} \in E^1$ s.t. $o(f_{n+1}) = z_n$, $z_{n+1} :=
t(f_{n+1})$ gets to $\lambda$ infinitely often, and $z_{n+1}
\notin \lambda^0$.  Notice that the number of representatives of
$[\lambda]$ that begin with $t(e_1)$, by (2) above, is finite.
Therefore $t(e_1)$ does not get to any of the $z_i$'s, that is,
$t(e_1)$ does not reach the infinite path $f_1f_2 \ldots$. Thus
$E$ is not co-final. Therefore $C^*(E)$ is not simple.
\end{proof}

We are now ready to prove the first of the measure results.

\begin{theorem}\label{them.1} Let E be a row-finite directed
graph with no sinks. $C^*(E)$ is liminal iff $E$ satisfies
condition $(M)$.
\end{theorem}

\begin{proof}
Suppose $E$ satisfies $(M)$.  Let $I$ be a primitive ideal of
$C^*(E)$, let $H=\{v:p_v \in I\}$, and let $F = F(E \setminus H)$.
By Lemma \ref{lemm.1}, $F^0$ is a maximal tail, and [\cite{B}
Theorem 4.1 (b)] implies that $C^*(F) \cong C^*(E)/I_H$.

\noindent Case I. $I = I_H$.

Then $I_H$ is a primitive ideal, hence Lemma \ref{lemm.4} implies
that $F$ has no circuits. Using Lemma \ref{lemm.3}, let $\lambda$
be a hereditary infinite path s.t. $F^0 = \overline{\lambda^0}$.

$C^*(E)/I_H \cong C^*(F) = I _{\lambda} =
\overline{span}\{s_\alpha s_\beta^*: \alpha, \beta \in F^*,~~
\textrm{s.t.}~~ t(\alpha) = t(\beta) \in \lambda^0\}$.

By [\cite{K} Proposition 2.1] $I_\lambda$ is Morita equivalent to
$C^*(\lambda) \cong \mathcal{K}(\ell^2(\alpha^0))$.  Therefore
$C^*(E)$ is liminal.

\noindent Case II. $I_H \subsetneqq I$.

We will first prove that $F(E \setminus H)$ has a circuit. If $F =
F(E \setminus H)$ has no circuits then by Lemma \ref{lemm.3} $F^0
= \overline{\lambda^0}$ for some hereditary infinite path
$\lambda$. Therefore $C^*(F)$ is simple, implying that
$C^*(E)/I_H$ is simple.  But $I/I_H$ is a (proper) ideal of
$C^*(E)/I_H$ therefore $I/I_H = 0$ implying that $I = I_H$.  A
contradiction.

Hence $F$ must have a circuit, say $\alpha$. Lemma \ref{lemm.2}
implies that $F^0 = \overline{\alpha^0}$. Using [\cite{K}
Proposition 2.1], $C^*(F)$ is Morita equivalent to $C^*(\alpha)$
which is Morita equivalent to $C(\mathbb{T})$ which is liminal.
Therefore $C^*(E)/I_H$ is liminal.   Since $I/I_H$ is a primitive
ideal of $C^*(E)/I_H$ we get $C^*(E)/I ~~ \cong ~~ C^*(E)/I_H ~~
\big/ ~~ I/I_H ~~ \cong ~~ \mathcal{K}$.  Hence $C^*(E)$ is
liminal.

To prove the converse, suppose $E$ does not satisfy condition
$(M)$, i.e., there exist an infinite path $\lambda$ and a
$v_\lambda \in E^0$ s.t. the number of representatives of
$[\lambda]$ that begin with $v_\lambda$ is infinite.

Suppose that $E$ has a non-terminal circuit, say $\alpha$. Let $v$
be a vertex of $\alpha$ s.t. $\exists e \in E^1$ which is not an
edge of $\alpha$ and $o(e) = v$. $p_v = s_\alpha^*s_\alpha \sim
s_\alpha s_\alpha^* < s_\alpha s_\alpha^* + s_es_e^* \leq p_v$.
Therefore $p_v$ is an infinite projection. Hence $C^*(E)$ can not
be liminal.

Suppose now that all circuits of $E$ are terminal and that the
number of representatives of $[\lambda]$ that begin with
$v_\lambda$ is infinite.  We might assume that $v_\lambda$ =
$o(\lambda)$.  We want to prove that $C^*(E)$ is not liminal.  If
$v$ is a vertex s.t. $V(v)$ does not intersect $\lambda^0$, we can
factor $C^*(E)$ by the ideal generated by $\{v\}$. Hence we might
assume that $\forall v \in E^0~~ v \geq \lambda^0$.  Moreover,
this process gets rid of any terminal circuits, and hence we may
assume that $E$ has no circuits.

Also, since $V(v_\lambda)$ is hereditary, by [\cite{K} Proposition
2.1], $I_{V(v_\lambda)}$ is Morita equivalent to
$C^*(E(v_\lambda))$. Therefore it suffices to show that
$C^*(E(v_\lambda))$ is not liminal. Hence we might assume that $E
= E(v_\lambda)$.

If $\forall v \in \lambda^0~~ \exists w \in \lambda^0~~
\textrm{s.t.}~~ v \geq w~~ \textrm{and}~~ |\{e \in E^1:o(e)~~ =~~
w\}| \geq 2$, then by Lemma \ref{lemm.5} $C^*(E)$ is not type I,
therefore it is not liminal.

Suppose $\exists u \in \lambda^0 ~~\textrm{s.t.}~~ \forall w \in
\lambda^0~~ \textrm{with}~~ u \geq w, |\{e \in E^1:o(e) = w\}| =
1.$ Notice that there is exactly one representative of $[\lambda]$
that begins with $u$.

By re-selecting $v_\lambda$ further along on $\lambda$, we might
assume that $\forall w \in \lambda^0 \setminus \{v_\lambda\}$ the
number of representatives of $[\lambda]$ that begin with $w$ is
finite.

Thus $E$ satisfies the following conditions:

\begin{enumerate}
\item $E$ has no circuits.

\item The number of representatives of $[\lambda]$ that begin with
$v_\lambda$ is infinite.

\item $v_\lambda$ is the only such vertex in $\lambda^0$.

\item $E = E(v_\lambda)$, and

\item  $\forall v \in E^0~~ \exists w \in \lambda^0~~
\textrm{s.t.}~~ v \geq w~~ (i.e., E^0 \geq \lambda^0)$.
\end{enumerate}

\noindent Therefore by Lemma \ref{lemm.51} we get:

i) $\{0\}$ is a primitive ideal of $C^*(E)$.

ii) $C^*(E)$ is not simple.

If $C^*(E)$ is liminal, by i), since $\{0\}$ is a primitive ideal
of $C^*(E)$, $C^*(E) \cong C^*(E)/\{0\}$ is *-isomorphic to
$\mathcal{K}$. But from ii) $C^*(E)$ can not be *-isomorphic to
$\mathcal{K}$ because $\mathcal{K}$ is a simple $C^*$-algebra.
Therefore $C^*(E)$ can not be liminal.  This concludes the proof
of the theorem.
\end{proof}

\section{The largest Liminal Ideal of $C^*$-Algebras of Graphs
with no singular vertices.}

In this section we will investigate a method of extracting the
largest liminal ideal of the $C^*$-algebra of a row finite graph
$E$ with no sinks.

Before we state the proposition, we will extend the definition of
the equivalence $\sim$ from $E^\infty$ to $E^{**} = E^\infty \cup
E^*$, as it is done in [\cite{S} Remark 1.10].  For $p, q \in
E^*$, we say $p \sim q$ if $t(p) = t(q)$.

The proposition gives a method of extracting the largest liminal
ideal of $C^*(E)$ of a graph $E$ with no singular vertices by
giving a characterization on the set of vertices that generate the
ideal. The first part of the proposition, which will eventually be
needed, can be proven for a general graph without much
complication. Therefore we state that part of the proposition for
a general graph.

\begin{proposition}\label{prop.1}  Let $E$ be a directed graph
and $H = \{v \in E^0: \forall [\lambda] \in (E^\infty \cup
\Lambda^*)/\sim$, the number of representatives of $[\lambda]$
that begin with $v$ is finite\}. Then

\begin{enumerate}

\item[(a)] $H$ is hereditary and saturated.

\item[(b)] If $E$ is row-finite with no sinks then $I_H$ is the
largest liminal ideal of $C^*(E)$.
\end{enumerate}
\end{proposition}

\begin{proof}
Suppose $v \in H$ and $v \geq w$.  Let $p$ be a path from $v$ to
$w$ and let $\lambda \in E^\infty \cup \Lambda^*$.  If $\beta \sim
\lambda$ and $o(\beta) = w$ then $p\beta \sim \lambda$ and
$o(p\beta) = v$. Therefore the number of representatives of
$[\lambda]$ that begin with $w$ is less than or equal to the
number of representatives of $[\lambda]$ that begin with $v$.
Therefore $w \in H$. Thus $H$ is hereditary.

Suppose $v \in E^0$ is not singular and $\{w \in E^0: v \geq w\}
\subseteq H$. Let $\bigtriangleup(v) = \{e \in E^1:o(e) = v\}$.
Note that $\bigtriangleup(v)$ is a finite set and $\forall e \in
\bigtriangleup(v),~~ t(e) \in H$. Let $\lambda \in E^\infty \cup
\Lambda^*$ and $\beta \sim \lambda$ where $o(\beta) = v$.  Then
the first edge of $\beta$ is in $\bigtriangleup(v)$.  Therefore
the number of representatives of $[\lambda]$ that begin with $v$
is equal to the sum of the number of representatives of
$[\lambda]$ that begin with a vertex in $\{t(e):e \in
\bigtriangleup(v)\}$, which is a finite sum of finite numbers.
Therefore $v \in H$. Hence $H$ is saturated.

To prove (b), suppose $E$ is row-finite with no sinks.  Let $F =
Graph(H)$.  Clearly $F$ satisfies condition $(M)$. Hence Theorem
\ref{them.1} implies that $C^*(F)$ is liminal.  By [\cite{K}
Proposition 2.1], $I_H$ is Morita equivalent to $C^*(F)$.  Hence
$I_H$ is a liminal ideal.  What remains is to prove that $I_H$ is
the largest liminal ideal of $C^*(E)$.

Let $I$ be the largest liminal ideal of $C^*(E)$. Thus $I_H
\subseteq I$. Since the largest liminal ideal of a $C^*$-algebra
is invariant under automorphisms, $I$ is gauge invariant,
therefore $I = I_K$ for some saturated hereditary subset $K$ of
$E^0$ which includes $H$. We will prove that $K \subseteq H$. Let
$G = Graph(K)$. Since $I_K$ is Morita equivalent to $C^*(G)$,
$C^*(G)$ is liminal hence, by Theorem \ref{them.1}, $G$ satisfies
condition $(M)$. Let $v \in K = G^0$.  If $\beta \in E^\infty$
with $o(\beta) = v$, because $K$ is hereditary, $\beta^0 \subseteq
K$. Therefore $\beta \in G^\infty$.  Now let $[\lambda] \in
E^\infty/\sim$,  and let $\gamma$ be a representative of
$[\lambda]$ that begins with $v$. (If no such $\gamma$ exists then
the number of representatives of $[\lambda]$ is zero.) Then
$\{\beta \in E^\infty: \beta \sim \lambda,~~ o(\beta) = v\} =
\{\beta \in G^\infty: \beta \sim \gamma,~~ o(\beta) = v\}$, i.e.,
the set of representatives of $[\lambda]$ that begin with $v$ is
subset of the set of representatives of $[\gamma]$ (as an
equivalence class of $G^\infty/\sim$) that begin with $v$. Since
$G$ satisfies condition $(M)$ the second set is finite. Therefore
$v \in H$, implying that $K \subseteq H$. Therefore $I_H = I_K$.
\end{proof}

\section{Liminal $C^*$-Algebras of General Graphs.}

In this section we will consider  for a general graph $E$ and give
the necessary and sufficient conditions for $C^*(E)$ to be liminal
in terms of the properties of the graph.

In \cite{Dr} the authors gave a recipe on how to ``desingularize a
graph $E$'', that is, obtain a graph $F$ that has no singular
vertices (by adding a tail at every singular vertex of $E$) so
that $C^*(E)$ and $C^*(F)$ are Morita equivalent.  Therefore, we
will use this idea of desingularizing $E$ and use the results of
the previous sections to get the needed results.

We will begin by extending the definition of condition $(M)$ from
row-finite graphs with no sinks to general graphs:

\begin{definition} \label{defn.2} A graph $E$ is said to satisfy
condition $(M)$ if $\forall [p] \in (E^\infty \cup \Lambda^*)
/\sim$ and any $v \in E^0$, the number of representatives of $[p]$
that begin with $v$ is finite.
\end{definition}

Notice that when $E$ is a row-finite graph with no sinks,
Definition \ref{defn.1} and Definition \ref{defn.2} say the same
thing.

Since we need to use the results of the previous sections, it is
important to check that condition $(M)$ is preserved by the
desingularization process.  We will do that in the next two
lemmas.

\begin{remark}\label{remk.2} [\cite{Dr} Lemma 2.6(a)] states that if
$F$ is a desingularization of a directed graph $E$ then there are
bijective maps.

$\phi: E^* \longrightarrow \{\beta \in F^*: o(\beta),~~ t(\beta)
\in E^0\}$

$\phi_\infty: E^\infty \cup \Lambda^* \longrightarrow \{\lambda
\in F^\infty: o(\lambda) \in E^0\}$

The map $\phi$ preserves origin and terminus (and hence preserves
circuits).  The map $\phi_\infty$ preserves origin.
\end{remark}

\begin{lemma}\label{lemm.7} The map $\phi_\infty$ preserves the
equivalence, in fact, for $p, q \in E^\infty \cup \Lambda^*$, $p
\sim q$ iff $\phi_\infty(p) \sim \phi_\infty(q)$.
\end{lemma}

\begin{proof}
Observe that if $\mu \nu \in E^\infty \cup \Lambda^*$ where $\mu
\in E^*$ then $\phi_\infty (\mu \nu) = \phi(\mu) \phi_\infty(\nu)$
and $\phi_\infty(\nu) \in F^\infty$

Now let  $p=e_1e_2 \ldots$, $q = f_1 f_2 \ldots \in E^\infty$ s.t.
$p \sim q,~~ \exists i, j$ s.t. $e_{i+r} = f_{j+r}~~ \forall r \in
\mathbb{N}$. Thus $p = \mu_1 \nu$ and $q = \mu_2 \nu$ where $\mu_1
= e_1 e_2 ... e_i,~~ \mu_2 = f_1 f_2... f_j$ and $\nu = e_{i+1}
e_{i+2} \ldots = f_{j+1} f_{j+2} \ldots$.  Therefore $\phi_\infty
(p) = \phi (\mu_1) \phi_\infty (\nu)$ and $\phi_\infty (q) =
\phi(\mu_2) \phi_\infty (\nu)$  implying $\phi_\infty(p) \sim
\phi_\infty(q)$.

If $p, q \in \Lambda^*$ s.t. $p \sim q$ then $t(p) = t(q)$  is a
singular vertex. Hence $\phi_\infty(t(p)) = \phi_\infty(t(q))$.
Moreover $\phi_\infty(p) = \phi(p)\phi_\infty(t(p))$ and
$\phi_\infty(q) = \phi(q)\phi_\infty(t(q))$ implying
$\phi_\infty(p) \sim \phi_\infty(q)$.

Hence $\phi_\infty(p) \sim  \phi_\infty(q)$ whenever $p \sim q$.

To prove the converse, suppose $\phi_\infty(p_1) \sim
\phi_\infty(p_2)$ for $p_1, p_2 \in E^\infty \cup \Lambda^*$.

Claim:  If $p_1 \in \Lambda^*$ then $p_2 \in \Lambda^*$. If $p_1
\in E^\infty$ then $p_2 \in E^\infty$.

We prove the claim.  Suppose $p_1 \in \Lambda^*$.  Thus
$\phi_\infty(p_1) = \phi(p_1)e_1e_2 \ldots$ where $e_1e_2 \ldots$
is the tail added to $t(p_1)$ in the construction of $F$, i.e.,
$t(p_1) = o(e_1e_2 \ldots)$. Therefore, $\phi_\infty(p_1) \sim
e_2e_3 \ldots$. Since $\phi_\infty(p_1) \sim \phi_\infty(p_2)$ we
get $\phi_\infty(p_2) \sim e_2e_3 \ldots$. If $p_2 \in E^\infty$
then $p_2 = f_1f_2 \ldots$ for some $f_1, f_2, \ldots \in E^1$
Therefore $\phi_\infty(p_2) = \phi(f_1) \phi(f_2) \ldots$.
Implying that $\phi(f_1) \phi(f_2) \ldots \sim e_2e_3 \ldots$. But
for each $i \geq 1$ we have $o(\phi(f_i)),~~ t(\phi(f_i)) \in E^0$ and by
the construction of $F$,  for each $i \geq 2$ $o(e_i), t(e_i)
\notin E^0$.  Therefore $\phi(f_1) \phi(f_2) \ldots$ can not be
equivalent to the path $e_2 e_3 \ldots$, which is a contradiction.
Therefore $p_2 \in \Lambda^*$.  The second statement follows from
the contrapositive of the first statement by symmetry.

Now suppose $p_1 \in \Lambda^*$.  By the above claim, $p_2 \in
\Lambda^*$.  Thus $\phi_\infty(p_2) = \phi(p_2)g_1g_2 \ldots$,
where $g_1g_2 \ldots$ is the tail added to $t(p_2)$ in the
construction of $F$.  Hence $\phi_\infty(p_2) \sim g_1g_2 \ldots$.
Since $\phi_\infty(p_1) \sim e_1e_2 \ldots$, we get $e_1e_2 \ldots
\sim g_1g_2 \ldots$. Notice that (by the construction of $F$)
$t(p_1)$ is the only entrance of $e_1e_2 \ldots$ and $t(p_2)$ is
the only entrance to $g_1g_2 \ldots$. Therefore either $t(p_1) =
o(g_i)$ for some $i$ or $t(p_2) = o(e_i)$ for some $i$.  WLOG
suppose $t(p_1) = o(g_i)$, thus $e_1e_2 \ldots = g_ig_{i+1}
\ldots$. But $t(p_2) = o(g_1)$ is the only vertex in the path
$g_1g_2 \ldots$ that belongs to $E^0$ and $t(p_1) \in E^0$ Hence
$t(p_1) = t(p_2)$.  Therefore $p_1 \sim p_2$.

If $p_1 \in E^\infty$ then, by the above claim, $p_2 \in
E^\infty$. Notice that $\forall v \in \phi_\infty(p_i)^0$ either
$v \in E^0$ (hence in $p_i^0$) or $\exists w \in p_i^0$ s.t. $v
\geq w$. Since $\phi_\infty(p_1) \sim \phi_\infty(p_2),~~
\phi_\infty(p_i) = \mu_i \nu$, for some $\mu_1, \mu_2 \in F^*$ and
some $\nu \in F^\infty$, and $t(\mu_1) = t(\mu_2) = o(\nu)$.
Extending $\mu_1$ and $\mu_2$ along $\nu$, if needed, we may
assume that $t(\mu_i) \in E^0$, i.e., $\mu_1,~~ \mu_2 \in \{\beta
\in F^*:o(\beta), t(\beta) \in E^0\},~~ \nu \in \{\beta \in
F^\infty:o(\beta) \in E^0\}$ and $t(\mu_1) = t(\mu_2) = o(\nu)$.
Therefore $\mu_i = \phi(\delta_i),~~ \nu = \phi_\infty(\gamma)$
for some $\delta_i \in E^*$ and some $\gamma \in E^\infty \cup
\Lambda^*.$  Implying $p_i = \phi_\infty^{-1}(\mu_i \nu) =
\phi_\infty^{-1}(\phi(\delta_i)\phi_\infty(\gamma)) =
\phi_\infty^{-1}(\phi_\infty(\delta_i \gamma)) = \delta_i \gamma.$
Thus $p_1 \sim p_2$.
\end{proof}

\begin{lemma}\label{lemm.8}  Let $F$ be a disingularization of a
directed graph $E$. Then $E$ satisfies condition $(M)$ iff $F$
satisfies condition $(M)$.
\end{lemma}

\begin{proof}  We will prove the only if side.  Recall that $F$ has
no singular vertices.  Suppose $F$ does not satisfy condition
$(M)$. Let $v \in F^0$ and $[\lambda] \in F^\infty / \sim$ s.t.
the number of representatives of $[\lambda]$ that begin with $v$
is infinite.

If $v \notin E^0$ then $v$ is on an added tail to a singular
vertex $v_0$ of $E$ and there is (only one) path from $v_0$ to
$v$. Then the number of representatives of $[\lambda]$ that begin
with $v$ (in the graph $F$) is equal to the number of
representatives of $[\lambda]$ that begin with $v_0$ (in the graph
$F$). If the latter is finite then the first is finite, hence we
might assume that $v \in E^0$. Moreover, every path in $F^\infty$
is equivalent to one whose origin lies in $E^0$. Therefore we
might choose a representative $\lambda$ with $o(\lambda) \in E^0$.

The set of representatives of $[\lambda]$ that begin with $v$ is
$\{\beta \in F^\infty: o(\beta) = v~~ and~~ \lambda \sim \beta\}$.
Since $\phi_\infty$ is bijective, $\phi_\infty^{-1} \{\beta \in
F^\infty :o(\beta) = v~~ and~~ \lambda \sim \beta\}$ is an
infinite subset of $E^\infty \cup \Lambda^*$. As
$\phi_\infty^{-1}$ preserves origin and the equivalence,
$\phi_\infty^{-1}\{\beta \in F^\infty : o(\beta) = v~~ and~~ \beta
\sim \lambda\} = \{\phi_\infty^{-1}(\beta) \in E^\infty \cup
\Lambda^*: o(\phi_\infty^{-1}(\beta))= v~~ and~~
\phi_\infty^{-1}(\beta) \sim \phi_\infty^{-1}(\lambda)\}$. Thus
$[\phi_\infty^{-1}(\lambda)]$ has infinite representatives that
begin with $v$. Therefore $E$ does not satisfy condition $(M)$.

To prove the converse, suppose $E$ does not satisfy condition
$(M)$. Let $v \in E^0$ and $[p] \in (E^\infty \cup \Lambda^*) /
\sim$ s.t. the number of representatives of $[p]$ that begin with
$v$ is infinite.  The set of representatives of $[p]$ that begin
with $v$ is $\{q \in E^\infty \cup \Lambda^*: o(q) = v~~ and~~ q
\sim p\}$.  Since $\phi_\infty$ is bijective, $\phi_\infty \{q \in
E^\infty \cup \Lambda^*:o(q) = v~~ and~~ q \sim p\}$ is an
infinite subset of $F^\infty$. As $\phi_\infty$ preserves origin
and the equivalence, $\phi_\infty\{q \in E^\infty \cup \Lambda^* :
o(q) = v~~ and~~ q \sim p\} = \{\phi_\infty(q) \in F^\infty:
o(\phi_\infty(q)= v~~ and~~ \phi_\infty(q) \sim \phi_\infty(p)\}$.
Thus $[\phi_\infty(p)]$ has infinitely many representatives that
begin with $v$. Therefore $F$ does not satisfy condition $(M)$.
\end{proof}

We can now write the main theorem in its full generalities.

\begin{theorem}\label{them.3}  Let $E$ be a directed graph.
$C^*(E)$ is liminal iff $E$ satisfies condition $(M)$.
\end{theorem}

\begin{proof}
Let $F$ be a desingularization of $E$.  $E$ satisfies condition
$(M)$

\hspace{2.35in} iff $F$ satisfies condition $(M)$

\hspace{2.35in} iff $C^*(F)$ is liminal

\hspace{2.35in} iff $C^*(E)$ is liminal
\end{proof}

\section{The largest Liminal Ideal of $C^*$-Algebras of General Graphs.}

In this section we will identify the largest liminal ideal of
$C^*(E)$ for a general graph $E$.

We will, once again, follow the construction in \cite{Dr}.  For a
hereditary saturated subset $H$ of $E^0$, define:

$$B_H := \{v \in \Lambda:0 < |o^{-1}(v) \cap t^{-1}(E^0 \setminus
H)| < \infty\}.$$

Thus $B_H$ is the set of infinite emitters that point into $H$
infinitely often and out of $H$ at least once but finitely often.
In \cite{Dr} it is proven that the set $\{(H,S): H$ is a
hereditary saturated subset of $E^0$ and $S \subseteq B_H\}$ is a
lattice with the lattice structure $(H, S) \leq (H', S')$ iff $H
\subseteq H'$ and $S \subseteq H'\cup S'$.  Observe that, since
$B_H \cap H = \emptyset,~~ (H, S) \leq (H, S')$ iff $S \subseteq
S'$.

Let $E$ be a directed graph and $F$ be a disingularization of $E$,
let $H$ be a hereditary saturated subset of $E^0$, and let $S
\subseteq B_H$.  Following the construction in \cite{Dr}, define:

$$ \widetilde{H} := H \cup \{v_n \in F^0:v_n ~~ is~~ on~~ a~~
tail~~ added~~ to~~ a~~ vertex~~ in~~ H\}.$$

Thus $\widetilde{H}$ is the smallest hereditary saturated subset
of $F^0$ containing $H$.

Let $S \subseteq B_H$, and let $v_0 \in S$.  Let $v_i = t(e_i)$,
where $e_1e_2 \ldots$ is the tail added to $v_0$ in the
construction of $F$.  If $N_{v_0}$ is the smallest non-negative
integer s.t. $t(e_j) \in H,~~ \forall j \geq N_{v_0}$,  we have
that $\forall j \geq N_{v_0}, ~~ v_j$ emits exactly two edges: one
pointing to $v_{j+1}$ and one pointing to a vertex in $H$.  Define

$$T_{v_0} := \{v_n \in F^0:v_n~~ is~~ on~~ a~~ tail~~ added~~ to~~
v_0~~  and~~  n \geq N_{v_0}\}$$ and $$\displaystyle{H_S :=
\widetilde{H} \cup \bigcup_{v_0 \in S}T_{v_0}}.$$

[\cite{Dr} Lemma 3.2] states that the above construction defines a
lattice isomorphism from the lattice $\{(H,S):H ~~ is~~ a~~
hereditary~~ saturated~~ subset~~ of~~ E^0~~ and~~ S \subseteq
B_H\}$ onto the lattice of hereditary saturated subsets of $F^0$.

Let $\{t_e, q_v\}$ be a generating Cuntz-Krieger $F$-family and
$\{s_e, p_v\}$ be the canonical generating Cuntz-Krieger
$E$-family.  Let $p = \sum_{v \in E^0} q_v$.  Since $C^*(E)$ and
$C^*(F)$ are Morita equivalent via the imprimitivity bimodule
$pC^*(F)$ it follows that the Rieffel correspondence between
ideals in $C^*(F)$ and ideals in $C^*(E)$ is given by the map $I
\longmapsto pIp$.

Let $H$ be a hereditary saturated subset of $E^0$ and $S \subseteq
B_H$.  For $v_0 \in S$, define $$\displaystyle{p_{v_0}^H :=
p_{v_0} - \mathop{\sum_{o(e) = v_0}}_{t(e) \notin H} s_es_e^*}$$
and $$I_{(H,S)} :=~~ the~~ ideal~~ generated~~ by~~ \{p_v : v \in
H\} \cup \{p_v^H:v_0 \in S\}.$$  [\cite{Dr} Proposition 3.3]
states that if $E$ satisfies condition $(K)$:  every vertex of $E$
lies on either no circuits or at least two circuits, then
$pI_{H_S}p = I_{(H,S)}$.  The assumption that $E$ satisfies
condition $(K)$ was only used to make sure that all the ideals of
$F$ are gauge invariant.  Therefore whenever $I$ is a gauge
invariant ideal of $C^*(F)$ and $H_S = \{v \in F^0:p_v \in I\}$,
since $I =I_{H_S}$, we have $pIp = pI_{H_S}p = I_{(H,S)}$.
Moreover [\cite{Dr} Theorem 3.5] states that if $E$ satisfies
condition $(K)$ then the map $(H,S) \longmapsto I_{(H,S)}$ is a
bijection from the lattice $\{(H,S):H~~ is ~~ a ~~ hereditary ~~
saturated ~~ subset ~~ of ~~ E^0 ~~ and ~~ S \subseteq H\}$ onto
the lattice of ideals in $C^*(E)$.  Without the assumption that
$E$ satisfies condition $(K)$ the bijection will be from the
lattice $\{(H,S):H~~ is ~~ a ~~ hereditary ~~ saturated ~~ subset
~~ of ~~ E^0 ~~ and ~~ S \subseteq H\}$ onto the lattice of gauge
invariant ideals in $C^*(E)$.  Hence the gauge invariant ideals of
$E$ are of the form $I_{(H,S)}$ for some hereditary saturated
subset $H$ of $E^0$ and for some $S \subseteq B_H$.

To identify the largest liminal ideal of $C^*(E)$, first recall
that the largest liminal ideal of a $C^*$-algebra is invariant
under automorphisms.  Therefore the largest liminal ideal of
$C^*(E)$ has to be of the form $I_{(H, S)}$ for some hereditary
saturated subset $H$ of $E^0$ and a subset $S$ of $B_H$.  We set
$H_l = \{v \in E^0: \forall [\lambda] \in (E^\infty \cup
\Lambda_E^*)/\sim,~~ the~~ number~~ of~~ representatives~~ of~~
[\lambda]~~ that~~ begin~~ with~~ v~~ is~~ finite\}$. Since
$Graph(H_l)$ satisfies condition $(M)$, we see that the ideal
$I_{H_l} = I_{(H_l, \emptyset)}$ is a subset of the largest
liminal ideal of $C^*(E)$. While it is true that $H = H_l$, as
illustrated in the following example, it is not automatically
clear what $S$ can be.

\begin{eg}
Consider the following graphs:
 $$
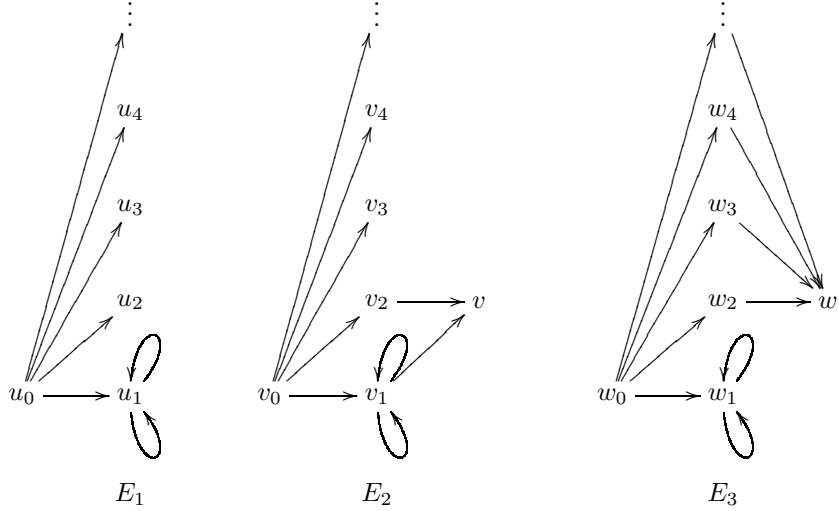
\xymatrix{  & {\vdots}\\
              & u_4 \\
              & u_3 \\
              & u_2\\
            u_0 \ar[ruuuu] \ar[ruuu] \ar[ruu] \ar[ru] \ar[r] &
                u_1 \ar@(d,dr)[] \ar@(ru,u)[]\\
                & E_1}
            \hspace{.5in}
 \xymatrix{  & {\vdots}\\
              & v_4 \\
              & v_3 \\
              & v_2 \ar[r] & v\\
            v_0 \ar[ruuuu] \ar[ruuu] \ar[ruu] \ar[ru] \ar[r] &
                v_1 \ar@(d,dr)[] \ar@(ru,u)[] \ar[ru]\\
                & E_2}
            \hspace{.5in}
\xymatrix{  & {\vdots} \ar[rddd] & \\
              & w_4 \ar[rdd] & \\
              & w_3 \ar[rd] & \\
              & w_2 \ar[r] & w\\
            w_0 \ar[ruuuu] \ar[ruuu] \ar[ruu] \ar[ru] \ar[r] &
                w_1 \ar@(d,dr)[] \ar@(ru,u)[]\\
                & E_3 & }$$

Let $I_{(H_i, S_i)}$ denote the largest liminal ideal of
$C^*(E_i)$. It is not hard to see that $H_1 = \{u_2, u_3, \ldots
\}$, $H_2 = \{v_2, v_3, \ldots \} \cup \{v\}$, $H_3 = \{w_2, w_3,
\ldots\} \cup \{w\}$, $B_{H_1} = \{u_0\}$, $B_{H_2} = \{v_0\}$,
and $B_{H_3} = \{w_0\}$.  A careful computation shows that $S_1 =
\{u_0\}$, $S_2 = \{v_0\}$ where as $S_3 = \emptyset$.  Notice that
we can reach from $v_0$ to $v$ in an infinite number of ways, but
not through $H$.  We can reach from $w_0$ to $w$ through $H$ in an
infinite number of ways.
\end{eg}

For a hereditary and saturated subset $H$ of $E^0$ and $v \in
B_H$, we define $D_{(v, H)} := \{e \in \Delta(v):t(e) \notin H\}$,
that is, $D_{(v,H)}$ is the set of all edges that begin with $v$
and point outside of $H$. Notice that $D_{(v, H)}$ is a non empty
finite set.

\begin{proposition}\label{prop.2}  Let $E$ be a directed graph
and $H = \{v \in E^0: \forall [\lambda] \in (E^\infty \cup
\Lambda_E^*)/\sim$, the number of representatives of $[\lambda]$
that begin with $v$ is finite\}.  Let $S = \{v \in B_H: E(v;
D_{(v,H)})$ satisfies condition $(M)\}$. Then $I_{(H,S)}$ is the
largest liminal ideal of $C^*(E)$.
\end{proposition}
\begin{proof}
That $H$ is hereditary and saturated is proven in Proposition
\ref{prop.1}.  Let $I_{(H',S')}$ be the largest liminal ideal of
$C^*(E)$ and let $F$ be a desingularization of $E$.  In what
follows, we will prove that $I_{(H,S)}$ = $I_{(H',S')}$.  To do
that we will prove: $H \subseteq H'$, $H' \subseteq H$, $S
\subseteq S'$ and $S' \subseteq S$, in that order.

We will prove that $H \subseteq H'$.   Notice that $I_{{H'}_{S'}}$
is the largest liminal ideal of $C^*(F)$. Using Proposition
\ref{prop.1} we get that ${H'}_{S'} = \{v \in F^0: \forall
[\lambda] \in F^\infty /\sim,$ the number of representatives  of
$[\lambda]$ that begin with $v$ is finite\}.

Let $G_H = Graph(H)$.  Notice that by [\cite{K} Proposition 2.1]
$I_H$ is Morita equivalent to $C^*(G_H)$.  Since $G_H$ satisfies
condition $(M)$, by Theorem \ref{them.1}, $C^*(G_H)$ is liminal.
Therefore $I_H = I_{(H,\emptyset)}$ is liminal. By the maximality
of $I_{(H', S')},~~ I_{(H,\emptyset)} \subseteq I_{(H', S')}$,
implying that $H \subseteq H'$.

We will prove that $H' \subseteq H$.  Let $G_{H'} = Graph(H')$.
$I_{H'} = I_{(H',\emptyset)} \subseteq I_{(H', S')}$.  Hence
$I_{H'}$ is liminal.  And by [\cite{K} Proposition 2.1], $I_{H'}$
is Morita equivalent to $C^*(G_{H'})$. Therefore $G_{H'}$
satisfies condition $(M)$.

Let $v \in H'$.  If $\beta \in E^{**}$ with $o(\beta) = v$ then,
since $H'$ is hereditary, $\beta \in G_{H'}$.  Now let $[\lambda]
\in (E^\infty \cup \Lambda^*)/\sim.$  If $\gamma$ is a
representative of $[\lambda]$ with $o(\gamma) = v$ then $\gamma
\in G_{H'}^\infty \cup \Lambda_{G_{H'}}^*$.  Therefore the set of
representatives of $[\lambda]$ that begin with $v$ is $\{\beta \in
E^\infty \cup \Lambda^*:o(\beta) = v, \beta \sim \gamma\} =
\{\beta \in G_{H'}^\infty \cup \Lambda_{G_{H'}}^*:o(\beta) = v,
\beta \sim \gamma\}$ which is finite, since $G_{H'}$ satisfies
condition $(M)$.  Therefore $v \in H$, hence $H' \subseteq H$.

Next we will prove that $S \subseteq S'$.  Let $v_0 \in S$.  To
show that $v_0 \in S'$ we will show that $v_n \in H_{S'}$ whenever
$n \geq N_{v_0}$, i.e., $\forall n \geq N_{v_0}$, and $\forall
[\lambda] \in F^\infty /\sim$, the number of representatives of
$[\lambda]$ that begin with $v_n$ is finite.

Let $n \geq N_{v_0}$ and let $[\lambda] \in F^\infty/\sim$.  If
$[\lambda]$  has no representative that begins with $v_n$ then
there is nothing to prove.  Let $\gamma$ be a representative of
$[\lambda]$ whith $o(\gamma) = v_n$.

First suppose that $\gamma^0 = \{v_n, v_{n+1}, \ldots\}$, i.e.,
$\gamma$ is the part of the tail added to $v_0$ in the
construction of $F$. Then $\{\beta \in F^\infty:o(\beta) = v_n,
\beta \sim \gamma\} = \{\gamma\}$ since $\gamma$ has no entry
other than $v_n$.  Therefore the number of representatives of
$[\lambda]$ that begin with $v_n$ is 1.

Now suppose $\gamma^0$ contains a vertex not in $\{v_n, v_{n+1},
\ldots\}$.  Recalling that $\forall k \geq N_{v_0},~~ v_k$ emits
exactly two edges, one pointing to $v_{k+1}$ and one pointing to a
vertex in $H$, let $w \in H$ be the first such vertex, i.e., $w
\in H \cap \gamma^0$ is chosen so that whenever $v \geq w$ and $v
\in \{v_n, v_{n+1}, \ldots\}$ then $v \notin H$. If $p$ is the
(only) path from $v_0$ to $v_n$ and $q$ is the path from $v_n$ to
$w$ along $\gamma$, then $\gamma = q\mu$ for some $\mu \in
F^\infty$ with $o(\mu) = w$. Moreover, $\phi_\infty^{-1}(p\gamma)
= \phi_\infty^{-1}(pq\mu) = \phi^{-1}(pq)\phi_\infty^{-1}(\mu)$
and $\phi^{-1}(pq)$ is an edge in $E^1$ with $o(\phi^{-1}(pq)) =
v_0$ and $t(\phi^{-1}(pq)) = w \in H$.  Therefore
$\phi_\infty^{-1}(p\gamma) \in E(v_0; D_{(v_0,H)})^\infty \cup
\Lambda_{E(v_0; D_{(v_0,H)})}^*$.  The set of representatives of
$[\lambda]$ that begin with $v_n$ is $\{\beta \in F^\infty :
o(\beta) = v_n, \beta \sim p\gamma\}$.  If $\beta \in F^\infty$ is
any representative of $[\lambda]$ that begins with $v_n$ then
$\beta \sim p\gamma \sim \mu$.  Hence $\beta^0$ has to contain a
vertex in $H$.  Applying the same argument on $\beta$ we see that
$p\beta$ is a representative of $[\lambda],~~ o(p\beta) = v$ and
$\phi_\infty^{-1}(p\beta) \in E(v_0; D_{(v_0,H)})^\infty \cup
\Lambda_{E(v_0; D_{(v_0,H)})}^*$.

Hence $|\{\beta \in F^\infty : o(\beta) = v_n, \beta \sim
p\gamma\}| = |\{p\beta \in F^\infty : p\beta \sim p\gamma\}| =
|\{\phi_\infty^{-1}(p\beta) \in E(v_0; D_{(v_0,H)})^\infty \cup
\Lambda^*_{E(v_0; D_{(v_0,H)})} : \phi_\infty^{-1}(p\beta) \sim
\phi_\infty^{-1}(p\gamma)\}|$ which is finite, since $E(v_0;
D_{(v_0,H)})$ satisfies condition $(M)$.

In each case, the number of representatives of $[\lambda]$ that
begin with $v_n$ is finite, implying that $v_n \in H_{S'}$.
Therefore $v_0 \in S'$.

Finally we will prove that $S' \subseteq S$.  Let $v_0 \in S'$.
We will show that $E(v_0; D_{(v_0,H)})$ satisfies condition $(M)$.
Let $\lambda \in E(v_0; D_{(v_0,H)})^\infty \cup \Lambda^*_{E(v_0;
D_{(v_0,H)})}$.  If a vertex $v \neq v_0$ is in $E(v_0;
D_{(v_0,H)})^0$ then it is in $H$, hence, by the definition of
$H$, the number of representatives of $[\lambda]$ that begin with
$v$ is finite.  What remains is to show that the number of
representatives of $[\lambda]$ that begin with $v_0$ is finite.
Noting that $v_{N_{v_0}} \in H_{S'}$, for any $\gamma \in
F^\infty$ the set $\{\mu \in F^\infty : o(\mu) = v_{N_{v_0}}, \mu
\sim \gamma\}$ is finite.  In particular, the set $\{\mu \in
F^\infty : o(\mu) = v_{N_{v_0}}, \mu \sim \phi_\infty(\lambda)\}$
is finite.

Let $\beta = e_1e_2 \ldots \in E(v_0; D_{(v_0,H)})^\infty \cup
\Lambda^*_{E(v_0; D_{(v_0,H)})}$ with $o(\beta) = v_0$.  Then
$\phi_\infty(\beta) = \phi(e_1)\phi_\infty(e_2e_3\ldots) \in
F^\infty$ and $o(\phi(e_1)) = v_0,~~ t(\phi(e_1)) =
o(\phi_\infty(e_2e_3\ldots)) \in H$.  Let $p$ be the path from
$v_0$ to $v_{N_{v_0}}$.

We will first show that the set $\{\beta = e_1e_2 \ldots \in
E(v_0; D_{(v_0,H)})^\infty \cup \Lambda^*_{E(v_0;
D_{(v_0,H)})}:e_1e_2\ldots \sim \lambda~~ and~~ v_{N_{v_0}} \in
\phi(e_1)^0\}$ is finite.

If $v_{N_{v_0}} \in \phi(e_1)^0$ then $\phi_\infty(\beta) = p\mu$
for some $\mu \in F^\infty$ with $o(\mu) = v_{N_{v_0}}$.  Hence
$|\{e_1e_2 \ldots \in E(v_0; D_{(v_0,H)})^\infty \cup
\Lambda^*_{E(v_0; D_{(v_0,H)})}:e_1e_2\ldots \sim \lambda~~ and~~
v_{N_{v_0}} \in \phi(e_1)^0\}| = |\{\phi_\infty(e_1e_2\ldots) \in
F^\infty:\phi_\infty(e_1e_2\ldots) \sim \phi_\infty(\lambda),
o(e_1) = v_0, t(e_1) \in H~~ and~~ v_{N_{v_0}} \in \phi(e_1)^0\}|
= |\{p\mu \in F^\infty: p\mu \sim \phi_\infty(\lambda)\}| = |\{\mu
\in F^\infty: o(\mu) = v_{N_{v_0}}~~ and~~ \mu \sim
\phi_\infty(\lambda)\}|$ which is finite.

We will next show that the set $\{e_1e_2 \ldots \in E(v_0;
D_{(v_0,H)})^\infty \cup \Lambda^*_{E(v_0;
D_{(v_0,H)})}:e_1e_2\ldots \sim \lambda~~ and~~ v_{N_{v_0}} \notin
\phi(e_1)^0\}$ is finite.

Observe that the set $\mathcal{E}:=\{e \in \Delta:t(e) \in H~~
and~~ v_{N_{v_0}} \notin \phi(e)\}$ is finite.  And $\forall e \in
\mathcal{E}$ the set $\{\beta \in E^\infty \cup \Lambda^*_E :
o(\beta) = t(e), \beta \sim \lambda\}$ is finite, since $\{t(e):e
\in  \mathcal{E}\} \subseteq H$.  Hence $|\{e_1e_2 \ldots \in
E(v_0; D_{(v_0,H)})^\infty \cup \Lambda^*_{E(v_0;
D_{(v_0,H)})}:e_1e_2\ldots \sim \lambda, v_{N_{v_0}} \notin
\phi(e_1)^0\}| = |\{e_1e_2 \ldots \in E^\infty \cup
\Lambda^*_E:e_1e_2\ldots \sim \lambda, t(e_1) \in K\}| = |\{\beta
\in E^\infty \cap \Lambda_E^*:o(\beta) \in K, \beta \sim
\lambda\}|$ which is finite, as the set is a finite union of
finite sets.

Therefore the set $\{\beta \in E(v_0; D_{(v_0,H)})^\infty \cup
\Lambda^*_{E(v_0; D_{(v_0,H)})}: \beta \sim \lambda\}$ is a union
of two finite sets, hence is finite. Thus $v_0 \in S$. It follows
that $S \subseteq S'$ concluding the proof.
\end{proof}

\section{Type I graph $C^*$-Algebras.}
In this section we will characterize Type I graph $C^*$-algebras.

We say that an edge $e$ reaches a path $p$ if t(e) reaches p, i.e.
if there is a path $q$ s.t. $o(q) = t(e)$ and $q \sim p$.

If $v$ is a sink then we regard $\{v\}$ as a tree.

For an infinite path $\lambda$, we use $\mathrm{N}_\lambda$ to
denote the number of vertices of $\lambda$ that emit multiple
edges that get back to $\lambda$.

\begin{lemma}\label{lemm.9} Let $E$ be a directed graph with:
\begin{enumerate}
\item Every circuit in $E$ is either terminal or transitory.
\item For any $\lambda \in E^\infty,~~ \mathrm{N}_\lambda$ is finite.
\end{enumerate}

Then $\exists v \in E^0$ s.t. $E(v)$ is either a terminal circuit
or a tree.
\end{lemma}

\begin{proof}
Let $z_1 \in E^0$.  If $E(z_1)$ is neither a terminal circuit nor
a tree, then there exists $z_2 \neq z_1$ s.t. $z_1$ and $z_2$ do
not belong to a common circuit, and there are (at least) two paths
from $z_1$ to $z_2$.

Notice that $\exists w_1 \in E^0$ s.t. $z_1 \geq w_1 \geq z_2$ and
$w_1$ emits multiple edges that reach $z_2$ (perhaps is $z_1$
itself). Observe that, by construction, $z_2 \ngeq z_1$.

Inductively: if $E(z_i)$ is neither a terminal circuit nor a tree,
then there exists $z_{i+1} \neq z_i$ s.t. $z_i$ and $z_{i+1}$ do
not belong to a common circuit, and there are (at least) two paths
from $z_i$ to $z_{i+1}$. Again $\exists w_i \in E^0$ s.t. $z_i
\geq w_i \geq z_{i+1}$ and $w_i$ emits multiple edges that reach
$z_{i+1}$. Observe also that $z_{i+1} \ngeq z_i$ and hence
$w_{i+1} \ngeq w_i$.

This process has to end, for otherwise, let $\lambda \in E^\infty$
be s.t. $\forall i~~ z_i,~~ w_i \in \lambda^0$. Then $\lambda$ has
infinite number of vertices that emit multiple edges that reach
$\lambda$, namely $w_1, w_2, \ldots$ contradicting the assumption.
\end{proof}

\begin{remark}\label{remk.3}  For $\lambda,~~ \gamma \in E^\infty$,
If $\lambda = p\gamma$, for some $p \in E^*$, then
$\mathrm{N}_\gamma \leq \mathrm{N}_\lambda \leq \mathrm{N}_\gamma
+ |p^0|$, where $|p^0|$ = the number of vertices in $p$, which is
finite since $p$ is a finite path. Therefore, $\mathrm{N}_\lambda$
is finite iff $\mathrm{N}_\gamma$ is finite. Moreover, if $\lambda
\sim \mu$ then $\lambda = p\gamma,~~ \mu = q\gamma$ for some $p,q
\in E^*$ and some $\gamma \in E^\infty$. Hence
$\mathrm{N}_\lambda$ is finite iff $\mathrm{N}_\gamma$ is finite
iff $\mathrm{N}_\mu$ is finite.
\end{remark}

\begin{theorem}\label{them.4}  Let $E$ be a graph.  $C^*(E)$ is
type I iff
\begin{enumerate}
\item Every circuit in $E$ is either terminal or transitory.
\item For any $\lambda \in E^\infty,~~ \mathrm{N}_\lambda$ is finite.
\end{enumerate}
\end{theorem}

We will first prove the following lemma.

\begin{lemma}\label{lemm.10} Let $E$ be a directed graph and $F$ be
a desingularization of $E$.  $E$ satisfies (1) and (2) of Theorem
\ref{them.4} iff $F$ satisfies (1) and (2) of Theorem \ref{them.4}
\end{lemma}
\begin{proof}  That $E$ satisfies (1) iff $F$ satisfies (1)
follows from the fact that the map $\phi$ of Remark \ref{remk.2}
preserves circuits.

Now we suppose that $E$ satisfies (1), equivalently $F$ satisfies
(1).

Suppose $E$ fails to satisfy (2).  Let $\lambda \in E^\infty$ s.t.
$\mathrm{N}_\lambda$ is infinite.  Suppose $v \in \lambda^0$ and
$p$ is a path s.t. $o(p) = v,~~ t(p) \in \lambda^0$.  Let $q$ be
the path along $\lambda$ s.t. $o(q) = v,~~ t(q) = t(p)$, then
$\exists \beta \in E^*$ and $\mu \in E^\infty$ s.t. $\lambda =
\beta q \mu$.  Since $\phi$ preserves origin and terminus,
$o(\phi(p)) = v = o(\phi(q))$ and $t(\phi(p)) = t(\phi(q))$.  And
$\phi_\infty(\lambda) = \phi_\infty(\beta q \mu) =
\phi(\beta)\phi(q)\phi_\infty(\mu).$  Since $\phi$ is bijective,
$\phi(p) = \phi(q)$ iff $p = q$.  Therefore, if $v$ (as a vertex
in $E$) emits multiple edges that get back to $\lambda$ then it
(as a vertex in $F$) emits multiple edges that get back to
$\phi(\lambda)$, implying that $\mathrm{N}_{\phi_\infty(\lambda)}$
is infinite.  Hence $F$ does not satisfy (2).

To prove the converse, suppose $E$ satisfies (2).  Let $\lambda
\in F^\infty$.  If $o(\lambda) \notin E^0$, then $o(\lambda)$ is
on a path extended from a singular vertex.  Using Remark
\ref{remk.3}, we may extend $\lambda$ (backwards) and assume that
$o(\lambda) \in E^0$.  Let $\gamma = {\phi_\infty}^{-1}(\lambda)
\in E^\infty \cup \Lambda^*$.

First suppose $\gamma \in \Lambda^*$.  Then $v_0:=t(\gamma) \sim
\gamma$. Hence $\phi_\infty(v_0) \sim \phi_\infty(\gamma) =
\lambda$. Using Remark \ref{remk.3}, we may assume that $\lambda =
\phi_\infty(v_0)$, that is, $\lambda$ is the path added to $v_o$
in the construction of $F$.  Thus each vertex of $\lambda$ emits
exactly two edges: one pointing to a vertex in $\lambda$ (the next
vertex) and one pointing to a vertex in $E^0$. Since $v_0$ is the
only entry to $\lambda$, if a vertex $v$ of $\lambda$ emits
multiple edges that get back to $\lambda$ then $v \geq v_0$.  And
since $F$ satisfies (1), there could be at most one such vector,
for otherewise $v_0$ would be on multiple circuits.  Hence
$\mathrm{N}_\lambda$ is at most 1.

Now suppose $\gamma \in E^\infty$.  Since $E$ satisfies (2),
$\mathrm{N}_\gamma$ is finite.  Going far enough on $\gamma$, let
$w \in \gamma^0$ be s.t. no vertex of $\gamma$ that $w$ can reach
to emits multiple edges that get back to $\gamma$. Let $\mu \in
E^*,~~ \beta \in E^\infty$ be s.t. $\gamma = \mu\beta$ and $t(\mu)
= w = o(\beta)$, then $\lambda = \phi(\mu)\phi_\infty(\beta)$.
Hence  $\lambda \sim \phi_\infty(\beta)$. Moreover, each $v \in
\beta^0$ emits exactly one edge that gets to $\beta$, which, in
fact, is an edge of $\beta$.

Let $v \in \beta^0$ and $p \in F^*$ be s.t. $o(p) = v, t(p) \in
\lambda^0$. Extending $p$, if needed, we may assume that $t(p) \in
\beta^0$. Let $q \in F^*$ be the path along $\lambda$ s.t. $o(q) =
v$ and $t(q) = t(p)$. Since $\phi$ is bijective, $\phi^{-1}(p) =
\phi^{-1}(q)$ iff $p = q$. But $v$ can get to $\beta$ in only one
way, therefore $\phi^{-1}(p) = \phi^{-1}(q)$, implying that $p =
q$.  Thus $v$ emits (in the graph $F$) only one edge that gets to
$\lambda$. Hence for each vertex $v \in \phi_\infty(\beta)$, if $v
\in E^0$ then $v$ emits only one edge that gets to $\lambda$.

Now let $v \in \phi_\infty(\beta) \setminus E^0$.  Then $v$ is on
a path extended from a singular vertex, say $v_0$.  Since $w \geq
v_0$, by the previous paragraph, $v_0$ emits only one edge that
gets to $\lambda$.  Let $p$ be the (only) path from $v_0$ to $v$.
Let $\mu, \nu \in F^*$ be s.t. $t(\mu), t(\nu) \in \lambda^0$ and
$o(\mu) = o(\nu) = v$ . Extending $\mu$ or $\nu$ along $\lambda$,
if needed, we can assume that $t(\mu) = t(\nu)$.  Again extending
them along $\lambda$ we can assume that  $t(\mu) = t(\nu) \in
\beta^0$. Observe that $o(p\mu) = o(p\nu) = v_0$ and $t(p\mu) =
t(p\nu) \in \beta^0$.  Therefore $o(\phi^{-1}(p\mu)) =
o(\phi^{-1}(p\nu)) = v_0$ and $t(\phi^{-1}(p\mu)) =
t(\phi^{-1}(p\nu)) \in \beta^0$. But each vertex in $\beta$ emits
exactly one edge that gets to $\beta$, i.e., there is exactly one
path from $v_0$ to $t(\phi^{-1}(p\mu))$ hence $p\mu = p\nu$.
Therefore, $\mu = \nu$.  That is,  $v$ emits only one edge that
gets to $\lambda$.  Therefore $\mathrm{N}_{\phi_\infty(\beta)} =
0$. By Remark \ref{remk.3} we get $\mathrm{N}_\lambda$ is finite.
\end{proof}

\begin{remark}
$E$ satisfies (2) of Theorem \ref{them.4} does not imply that its
desingularization $F$ satisfies (2) of Theorem \ref{them.4} as
illustrated by the following example.
\end{remark}

\begin{eg}
If $E$ is the $\mathcal{O}_\infty$ graph (one vertex with
infinitely many loops), which clearly satisfies (2) of Theorem
\ref{them.4}, then its disingularization does not satisfy (2) of
Theorem \ref{them.4}.  The disingularization looks like this:

$$
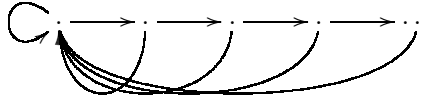
\xymatrix{.\ar@(ul,dl)[] \ar[r] & . \ar[r] \ar@(d,d)[l] & .
\ar[r] \ar@(d,d)[ll] & . \ar[r] \ar@(d,d)[lll] & {\ldots}
\ar@(d,d)[llll]}$$ \vspace{.1in}
\end{eg}

\begin{proof}[Proof of Theorem \ref{them.4}] We first prove the
if side. We will first assume that $E$ is a row-finite graph with
no sinks. Let $(I_\rho)_{0 \leq \rho \leq \alpha}$ be an
increasing family of ideals of $C^*(E)$ s.t.
\begin{enumerate}
\item[(a)]  $I_0 = \{0\},~~ C^*(E)/I_\alpha$ is antiliminal.
\item[(b)]  If $\rho \leq \alpha$ is a limit ordinal,
$\displaystyle{I_\rho = \overline{\bigcup_{\beta < \rho}I_\beta}}$
\item[(c)]  If $\rho < \alpha,~~ I_{\rho+1}/I_\rho$ is a liminal ideal
of $C^*(E)/I_\rho$ and is non zero.
\end{enumerate}

We prove that $I_\alpha = C^*(E)$.  Since $I_\alpha$ is the
largest Type I ideal of $C^*(E)$, it is gauge invariant. Let $H$
be a hereditary saturated subset of $E^0$ s.t. $I_\alpha = I_H$.
If $H \neq E^0$ then let $F = F(E \setminus H)$. Clearly $F$
satisfies (1) and (2) of the theorem. Using Lemma \ref{lemm.9} let
$v_0 \in F^0$ be s.t. $K=\{v \in F^0:v_0 \geq v\}$ is the set of
vertices of either a terminal circuit or a tree. Let $G =
Graph(K)$, thus $G$ is either a terminal circuit or a tree. By
[\cite{K} Proposition 2.1] $I_K$ is Morita equivalent to $C^*(G)$.
Moreover $G$ satifies condition $(M)$, hence by Theorem
\ref{them.1} $C^*(G)$ is liminal. And $I_K$ is an ideal of $C^*(F)
\cong C^*(E)/I_\alpha$ contradicting the assumption that
$C^*(E)/I_\alpha$ is antiliminal. It follows that $I_\alpha =
C^*(E)$. Therefore $C^*(E)$ is Type I.

For an arbitrary graph $E$, let $F$ be a desingularization of $E$.
By Lemma \ref{lemm.10} $F$ satisfies (1) and (2) of the theorem.
And by the above argument, $C^*(F)$ is Type I.  Therefore $C^*(E)$
is Type I.

To prove the converse, suppose $E$ has a non-terminal
non-transitory circuit, that is, $E$ has a vertex that is on (at
least) two circuits. Let $v_0$ be a vertex on two circuits, say
$\alpha$ and $\beta$. Let $F$ be the sub-graph containing (only)
the edges and vertices of $\alpha$ and $\beta$.

$\mathcal{A} := \overline{span}\{s_\mu s_\nu^*: \mu, \nu~~ are~~
paths~~ made~~ by~~ \alpha~~ and~~ \beta~~ or~~ just~~ v_0\}$  is
a $C^*$-subalgebra of $C^*(F)$.  But $\mathcal{A} \cong
\mathcal{O}_2$ which is not Type I.  Hence $C^*(F)$ is not Type I.
By Remark \ref{remk.1} $C^*(E)$ has a sub-algebra whose quotient
is not Type I therefore $C^*(E)$  is not Type I.

Suppose now that each circuit in $E$ is either terminal or
transitory and $\exists \lambda \in E^\infty$ s.t.
$\mathrm{N}_\lambda$ is infinite. Let $v_\lambda = o(\lambda)$.
Let $G = E(v)$.  If $v$ is a vertex s.t. $V(v)$ does not intersect
$\lambda^0$, we can factor $C^*(G)$ by the ideal generated by
$\{v\}$. This process gets rid of any terminal circuits of $G$. By
Lemma \ref{lemm.5} $C^*(G)$ is not Type I, implying that $C^*(E)$
is not Type I.
\end{proof}

Next we will identify the largest Type I ideal of the
$C^*$-algebra of a graph $E$.  For a vertex $v$ of $E$
(respectively $F$), recall that $E(v)$ (respectively $F(v)$)
denotes the sub-graph of $E$ (respectively $F$) that $v$ can
`see'.

We begin with the following lemma.

\begin{lemma}\label{lemm.11} Let $E$ be a directed graph, $F$ a
desingularization of $E$ and $v \in E^0$.  Then $F(v)$ is a
desingularization of $E(v)$.
\end{lemma}

\begin{proof}
Let $u \in E(v)^0 = \{w \in E: v \geq w\}$.  Let $p$ be a path in
$E$ with $o(p) = v$, and $t(p) = u$.  then $\phi(p)$ is a path in
$F$ with $o(\phi(p)) = v$, and $t(\phi(p)) = u$. Hence $u \in
F(v)^0$, implying that $E(v)^0 \subseteq F(v)^0$. Clearly $F(v)$
has no singular vertices.  Let $v_0 \in E(v)^0$ be a singular
vertex. If $v_n$ is a vertex on the path added to $v_0$ in the
construction of $F$, since $F(v)^0$ is hereditary and $v_0 \in
F(v)^0$, we get $v_n \in F(v)^0$.  Therefore the path added to
$v_0$ is in the graph $F(v)$.  To show that $F(v)$ has exactly the
vertices needed to desingularize $E(v)$, let $w \in F(v)^0$. Let
$p$ be a path in $F(v)$ with $o(p) = v$ and $t(p) = w$.  If $w \in
E^0$ then $\phi^{-1}(p) \in E^*$ and $o(\phi^{-1}(p)) = v$ and
$t(\phi^{-1}(p)) = w$.  Therefore $v \geq w$ in the graph $E$.
Hence $w \in E(v)^0$.  If $w \notin E^0$ then there is a singular
vertex, say $v_0 \in E^0$ s.t. $w$ is on the path added to $v_0$
in the construction of $F$.  Since the path from $v_0$ to $w$ has
no other entry than $v_0$ and since $v \geq w$, we must have $v
\geq v_0$.  Hence $w$ is on the the graph obtained when $E(v)$ is
desingularized.  Therefore $F(v)$ is a desingularization of
$E(v)$.
\end{proof}

The following corollary follows from Lemma \ref{lemm.11} and Lemma
\ref{lemm.10}.

\begin{corollary}\label{coro.1} Let $E$ be a directed graph, $F$ a
desingularization of $E$ and $v \in E^0$.  Then $E(v)$ satisfies
(1) and (2) of Theorem \ref{them.4} iff $F(v)$ satisfies (1) and
(2) of Theorem \ref{them.4}.

\end{corollary}

The next proposition identifies the largest Type I ideal of the
$C^*$-algebra of a row-finite graph $E$ with no sinks. The first
part of the proposition, which will be needed later, is written
for a general graph as it is proven without the need of the
property that $E$ is row-finite and has no sinks.

\begin{proposition}\label{prop.3}
Let $E$ be a directed graph and $$H = \{v \in E^0: E(v)~~
satisfies~~ (1)~~ and~~ (2)~~ of~~ Theorem~~ \ref{them.4} \}.$$
Then
\begin{enumerate}
\item[(a)] $H$ is a hereditary saturated subset of $E^0$.

\item[(b)]  If $E$ is a row-finite graph with no sinks then $I_H$
is the largest Type I ideal of $C^*(E)$.
\end{enumerate}
\end{proposition}
\begin{proof}
We first prove (a).  That $H$ is hereditary follows from $v \geq w
\Longrightarrow E(v) \supseteq E(w)$.  We prove that $H$ is
saturated.  Suppose $v \in E^0$ and $\{w \in E^0: v \geq w\}
\subseteq H$. Let $\bigtriangleup(v) = \{e \in E^1:o(e) = v\}$.
Note that $\forall e \in \bigtriangleup(v)$,  $t(e) \in H$. If
there is a circuit at $v$, i.e., $v$ is a vertex of some circuit,
then $v \geq v$, implying that $v \in H$. Suppose there are no
circuits at $v$.  If there is a vertex $w \in E(v)^0$ on a circuit
then it is in $E(t(e))^0$ for some $e \in \bigtriangleup(v)$. But
$t(e) \in H$, hence $w$ can not be on multiple circuits, i.e,
$E(v)$ has no non-terminal and non-transitory circuits. Hence
$E(v)$ satisfies (1) of Theorem \ref{them.4}.  Let $\lambda \in
E(v)^\infty$ then $\exists e \in \bigtriangleup(v)$ and $\beta \in
E(t(e))$ s.t. $\lambda \sim \beta$.  Since $t(e) \in H$,
$\mathrm{N}_\beta$ is finite.  Using Remark \ref{remk.3} we get
that $\mathrm{N}_\lambda$ is finite. Therefore $v \in H$.  Hence
$H$ is saturated.

To prove (b), suppose $E$ is row-finite with no sinks. Let $F =
Graph(H)$.  Clearly $F$ satisfies (1) and (2) of Theorem
\ref{them.4}, hence by Theorem \ref{them.4}, $C^*(F)$ is Type I.
Moreover, by [\cite{K} Proposition 2.1], $I_H$ is Morita
equivalent to $C^*(F)$. Hence $I_H$ is Type I. Let $I$ be the
largest Type I ideal of $C^*(E)$, then $I_H \subseteq I$.  Since
$I$ is gauge invariant, $I = I_K$ for some hereditary saturated
subset $K$ of $E^0$ that includes $H$. We will prove that $K
\subseteq H$.  Let $G = Graph(K)$.  Since $I_K$ is Morita
equivalent to $C^*(G)$, $C^*(G)$ is Type I, hence $G$ satisfies
(1) and (2) of Theorem \ref{them.4}.  Let $v \in K$, since $E(v)
\subseteq G$, $E(v)$ satisfies (1) and (2) of Theorem
\ref{them.4}. Therefore $v \in H$, hence $K \subseteq H$.
\end{proof}

The next proposition generalizes Proposition \ref{prop.3}.

\begin{proposition}\label{prop.4}
Let $E$ be a directed graph and $$H = \{v \in E^0: E(v)~~
satisfies~~ (1)~~ and~~ (2)~~ of~~ Theorem~~ \ref{them.4} \}.$$
Then $I_{(H,B_H)}$ is the largest Type I ideal of $C^*(E)$.

\end{proposition}

\begin{proof}

Let $I_{(H',S')}$ be the largest Type I ideal of $C^*(E)$ and let
$F$ be a desingularization of $E$ then $I_{{H'}_{S'}}$ is the
largest Type I ideal of $C^*(F)$. From (b) of Proposition
\ref{prop.3}, we get that ${H'}_{S'} = \{v \in F^0:F(v)~~
satisfies~~ (1)~~ and~~ (2)~~ of~~ Theorem~~ \ref{them.4}\}.$

We will prove that $H \subseteq H'$.  Let $G_H = Graph(H)$.
Clearly $G_H$ satisfies (1) and (2) of Theorem \ref{them.4} hence
$C^*(G_H)$ is Type I. By [\cite{K} Proposition 2.1] $I_H$ is
Morita equivalent to $C^*(G_H)$. Therefore $I_H =
I_{(H,\emptyset)}$ is Type I. By the maximality of $I_{(H',
S')},~~ I_{(H,\emptyset)} \subseteq I_{(H', S')}$, implying that
$H \subseteq H'$.

We will prove that $H' \subseteq H$.  Let $G_{H'} = Graph(H')$.
$I_{H'} = I_{(H',\emptyset)} \subseteq I_{(H', S')}$.  Hence
$I_{H'}$ is liminal. By [\cite{K} Proposition 2.1] $I_{H'}$ is
Morita equivalent to $C^*(G_{H'})$, implying that $C^*(G_{H'})$ is
liminal. Hence $G_{H'}$ satisfies (1) and (2) of Theorem
\ref{them.4}.

Let $v \in H'$.  Since $H'$ is hereditary and $E(v)^0 \subseteq
H'$ it follows that $E(v)$ is a sub-graph  of $G_{H'}$.  Thus
$E(v)$ satisfis (1) and (2) of Theorem \ref{them.4}.  Therefore $v
\in H$, hence $H' \subseteq H$.

Since $S' \subseteq B_H$, as $H ~~ = ~~ H'$, it remains to prove
that $B_H \subseteq S'$. Let $v_0 \in B_H$. To show that $v_0 \in
S'$ we will show that $ \forall n \geq N_{v_0}~~ v_n \in H_{S'}$
i.e., $F(v_n)$ satisfies (1) and (2) of Theorem \ref{them.4}. Let
$n \geq N_{v_0}$ and suppose $F(v_n)$ does not satisfy (1) of
Theorem \ref{them.4}.  Let $\alpha$ be a non-terminal and
non-transitory circuit in $F(v_n)$, and let $v \in \alpha^0$.

If $v$ is on the infinite path added to $v_0$ in the construction
of $F$ then $v_0$ is in the circuit $\alpha$. Notice that $v_n
\geq v \geq v_0$.  Recall that $\forall k \geq N_{v_0}~~ v_k$
emits exactly two edges one pointing to $v_{k+1}$ and one pointing
to a vertex in $H$.  Following along $\alpha$, we get that $v \geq
w$ for some vertex $w \in H$ of $\alpha$. But $H$ is hereditary,
therefore $v_0 \in H$, which contradicts to the fact that $H \cap
B_H = \emptyset$.

Suppose now that $v$ is not on the infinite path added to $v_0$.
Let $p$ be a path from $v_n$ to $v$. $p$ must contain a vertex,
say $w$, in $H$.  Notice that $w \geq v$ which implies that $v \in
F(w)$. Since $F(w)^0$ is hereditary, $\alpha$ is in the graph
$F(w)$. Hence $F(w)$ contains a non-terminal and non-transitory
circuit. Since $w \in H$, $E(w)$ satisfies (1) and (2) of Theorem
\ref{them.4}.  But this contradicts to Corollary \ref{coro.1}.
Therefore $F(v_n)$ satisfies (1) of Theorem \ref{them.4}.

To prove that $F(v_n)$ satisfies (2) of Theorem \ref{them.4}, let
$\lambda \in F(v_n)^\infty$.  Either $\lambda$ is on the tail
added to $v_0$ on the construction of $F$ or $\lambda^0$ contains
a vertex in $H$.

If $\lambda$ is on the tail added to $v_0$ then
$\mathrm{N}_\lambda = 0$.  Otherwise let $w \in \lambda^0 \cap H$.
Then $\lambda = p \mu$ for some $p \in F(v_n)^*$ and some $\mu \in
F(v_n)^\infty$ with $o(p) = v_n,~~ t(p) = w = o(\mu)$.  Implying
that $\lambda \sim \mu$. Since $w \in H$, $E(w)$ satisfies (1) and
(2) of Theorem \ref{them.4}.  By Corollary \ref{coro.1}, we get
that $F(w)$ satisfies (1) and (2) of Theorem \ref{them.4}.  Hence
$\mathrm{N}_\mu$ is finite and Remark \ref{remk.3} implies that
$\mathrm{N}_\lambda$ is finite.  Therefore $F(v_n)$ satisfies (2)
of Theorem \ref{them.4}.

We have established that $F(v_n)$ satisfies (1) and (2) of Theorem
\ref{them.4}. Therefore $v_n \in H_{S'}$ and hence $B_H \subseteq
S'$.  This concludes the proof.
\end{proof}

\bibliographystyle{amsplain}

\begin{thebibliography}{10}

\bibitem {B}  T. Bates, D. Pask, I. Reaburn and W. Szyma\'nski,
\textit{The $C^*$-Algebras of Row-Finite Graphs.} New York J.
Math. \textbf{6} (2000), 307--324.

\bibitem {C}  J. Cuntz, \textit{Simple $C^*$-Algebras Generated by
Isometries} Commun. Math. Phys. \textbf{57} (1977), 173--185.

\bibitem {Di}  J. Dixmier, \textit{$C^*$-Algebras.} North-Holland
Publishing Co., 1977.

\bibitem {Dr} D. Drinen and M. Tomforde, \textit{The
$C^*$-algebras of Arbitrary Graphs}, preprint(2001), Front for the
Mathematics ArXiv, math.OA/0009228.

\bibitem {K} A. Kumjian, D. Pask and I. Raeburn,
\textit{Cuntz-Krieger Algebras of Directed Graphs}, Pacific J.
Math. \textbf{184}, (1998), 161--174.

\bibitem {M}  G. Murphy, \textit{$C^*$-Algebras And Operator
Theory.} Academic Press, 1990.

\bibitem {Se}  J. Serre, \textit{Trees.} Springer-Verlag, Berlin, 1980.

\bibitem {S} J. Spielberg, \textit{A Functorial Approach to the
$C^*$-algebras of a Graph}, Internat. J. Math. \textbf{13}, (2002)
No.3, 245--277.

\end{thebibliography}

\end{document}